\documentclass{gtpart}   
\usepackage{pinlabel}
\usepackage[all]{xy}

\usepackage[utf8]{inputenc}
\usepackage{setspace}
\usepackage{graphicx}
\usepackage{graphics}
\usepackage[mathscr]{euscript}
\usepackage{tikz}
\usetikzlibrary{calc,trees,positioning,arrows,chains,shapes.geometric,
  decorations.pathreplacing,decorations.pathmorphing,shapes,matrix}
\usetikzlibrary{calc}
\usetikzlibrary{matrix}
\usepackage{mathtools}

\usepackage{enumitem}

\usepackage[all]{xy}
\usepackage{color}

\usepackage{bbold}

\usepackage{tikz}
\xyoption{all}
\usepackage{multirow}
\usepackage{etex, pictexwd,dcpic}
\usetikzlibrary{positioning}
\usetikzlibrary{shapes.geometric}
\usetikzlibrary{shapes.misc}
\usetikzlibrary{calc}
\usetikzlibrary{positioning}

\usetikzlibrary{trees} 
\usetikzlibrary[trees] 

\newtheorem{theorem}{Theorem}
\newtheorem{definition}[theorem]{Definition}
\newtheorem{lemma}[theorem]{Lemma}
\newtheorem{corollary}[theorem]{Corollary}
\newtheorem{proposition}[theorem]{Proposition}

\newtheorem{remark}{Remark}

\newtheorem{example}{Example}

\numberwithin{equation}{section}

\renewcommand{\(}{\begin{equation*}}
\renewcommand{\)}{\end{equation*}}
\newcommand{\bea}{\begin{eqnarray*}}
\newcommand{\eea}{\end{eqnarray*}}

\def\endofproof {\hfill{$\Box$}\\}

\newcommand{\cQ}{\ensuremath{\mathcal Q}}

\newcommand{\beq}{\begin{equation}}
\newcommand{\eeq}{\end{equation}}

\newcommand{\onto}{\twoheadrightarrow}

\newcommand{\into}{\hookrightarrow}

\newcommand{\theproof}{\noindent {\bf Proof.\ }}

\numberwithin{equation}{section}

\renewcommand{\(}{\begin{equation}}
\renewcommand{\)}{\end{equation}}

\def \qed{\mbox{ $\square$}}

\def\1{{\bf 1}}

\def\<{\langle}
\def\>{\rangle}

\numberwithin{equation}{section}


\newcommand{\RR}{\ensuremath{\mathbb R}}
\newcommand{\NN}{\ensuremath{\mathbb N}}

\newcommand{\ZZ}{\ensuremath{\mathbb Z}}
\newcommand{\QQ}{\ensuremath{\mathbb Q}}

\newcommand{\BB}{\ensuremath{\mathbf B}}

\newcommand{\chp}{\ensuremath{\mathscr{C}\mathrm{h}^{+}}}

\newcommand{\sset}{\ensuremath{s\mathscr{S}\mathrm{et}}}
\newcommand{\Top}{\ensuremath{\mathscr{T}\mathrm{op}}}

\newcommand{\ab}{\ensuremath{\mathscr{A}\mathrm{b}}}

\newcommand{\cartsp}{\mathscr{C}\mathrm{art}\mathscr{S}\mathrm{p}}

\newcommand{\E}{\ensuremath{\mathscr{E}}}

\newcommand{\map}{\mathrm{Map}}

\newcommand{\tK}{\widetilde{K}}
\newcommand{\tQ}{\widetilde{Q}}

\begin{document}

\title{Spectral sequences in smooth generalized cohomology}

\author{Daniel Grady}
\givenname{Daniel}
\surname{Grady}
\email{djg9@nyu.edu}
\address{Department of Mathematics\\
New York University\\\newline
Abu Dhabi\\
Abu Dhabi, UAE}

\author{Hisham Sati}
\givenname{Hisham}
\surname{Sati}
\email{hsati@nyu.edu}
\address{Department of Mathematics\\
New York University\\\newline
Abu Dhabi\\
Abu Dhabi, UAE}

\address{Department of Mathematics\\
University of Pittsburgh\newline
Pittsburgh, PA 15260, USA}

\keyword{Differential cohomology}
\keyword{Atiyah-Hirzebruch spectral sequence}
\keyword{Generalized cohomology} 
\subject{primary}{msc2010}{55S25}
\subject{secondary}{msc2010}{55N20}
\subject{secondary}{msc2010}{19L50}

\begin{abstract}
We consider spectral sequences in smooth generalized cohomology theories, including 
differential generalized cohomology theories. The main differential spectral sequences will be 
of the Atiyah-Hirzebruch (AHSS) type, where we provide a filtration by the {\v C}ech resolution
of smooth manifolds. This allows for systematic study of torsion in differential cohomology.
We apply this in detail to smooth Deligne cohomology, differential topological complex
 K-theory, and to a smooth extension of integral Morava K-theory that we introduce. In each 
 case we explicitly identify the differentials in the corresponding spectral sequences, which exhibit 
 an interesting and systematic interplay between (refinement of) classical cohomology operations,
 operations involving differential forms, and operations on cohomology with $U(1)$ coefficients. 
\end{abstract}

\maketitle

\tableofcontents

\section{Introduction} 

Spectral sequences are very useful algebraic tools that often allow for efficient computations
that would otherwise  require brute force (see \cite{Mc} for a broad survey). 
The Atiyah-Hirzebruch spectral sequence (henceforth AHSS) 
for K-theory and any generalized cohomology theory, in the topological sense,
was introduced  by Atiyah and Hirzebruch in \cite{AH}.  Excellent introduction to the generalized 
cohomology AHSS can also be found in Hilton \cite{Hi} and Adams \cite{Ad} (Sec. III.7). 
Other useful references on the subject include Switzer \cite{Sw} (Sec. 15, from a homology point of view, 
including the Gysin Sequence from AHSS) and interesting remarks in relation to spectra are given 
in Rudyak \cite{Ru} (Theorem 3.45 (homology), Remark 4.24 (Sheaves and Cech), Remark 
4.34 (Postnikov), and Corollary 7.12). A description with an eye for applications is 
given in \cite{HJJS} (Ch. 21).

The goal of this paper is to systematically study the spectral sequence in the context of 
smooth or differential cohomology (see \cite{CS} \cite{Fr} \cite{HS} \cite{SSu} 
\cite{Bun} \cite{BS} \cite{Urs}). Existence and interesting aspects of the AHSS in twisted forms of 
such differential cohomology theories have been considered briefly by Bunke and Nikolaus \cite{BN}, 
where the main interest was the effect of the geometric part of the twist on the spectral sequence. 
In this paper we take a step back 
and consider untwisted differential generalized cohomology to systematically study the 
corresponding AHSS 
in generality and determine  the differentials explicitly as cohomology operations. 
From the geometric point of view, one might expect on general grounds
that the geometric information carried by the differential cohomology theory should 
somehow manifest itself within the spectral sequence. On the other hand, from an algebraic 
point of view, one might a priori not expect much of that information to be retained, or 
expect it to even be totally stripped out while running through the homological algebra machine. 
We will show that the answer lies somewhat in between, and both intuitions are to some extent correct: 
The differentials in the spectral sequence will be essentially refinements of classical ones, but
with additional operations on differential forms. We recently characterized such operations  
 in \cite{GS2}, and so this paper is a natural continuation of that work.

Just as generalized cohomology theories are represented by spectra, differential cohomology
 theories are represented by certain sheaves of smooth spectra called \emph{differential function spectra}. 
 The original definition of differential function spectra was due to Hopkins and Singer in \cite{HS}, 
 generalized by Bunke, Nikolaus and V\"olkl in \cite{BNV} and was reformulated in terms of 
 cohesion by Schreiber in \cite{Urs}. The terms {\it smooth cohomology} and {\it differential cohomology} 
 seem to be used interchangeably in some of the literature (see e.g. \cite{BS2}). However, we will find it 
 useful for us to provide  a specific and precise usage, where the first is viewed as being more general 
 than the second. We also present most of our $\infty$-categories as combinatorial, simplicial model categories, rather than quasi-categories. We believe that this way nice objects are more easily and explicitly identifiable, which is 
 something desirable when dealing with differential cohomology. Indeed, our discussion will be very 
 explicit and the results will be readily utilizable.

 Ordinary cohomology has smooth extension with various different realizations, including those of 
\cite{CS} \cite{Ga} \cite{Bry} \cite{DL} \cite{HS} \cite{BKS}. All these realizations are in 
fact isomorphic \cite{SSu} \cite{BS2}. 
A description of K-theory with coefficients that combines vector bundles, connections, and differential 
forms into a topological context was initiated in \cite{Ka}. Using Karoubi's description Lott introduced 
$\R/\ZZ$-valued K-theory \cite{Lo} as well as   differential flat K-theory \cite{Lo2}. 
Currently there are various geometric models of differential K-theory \cite{Lo} \cite{BS} \cite{SSu} \cite{FL} \cite{TWZ1} \cite{TWZ2}. As in the case of ordinary differential cohomology these models should be equivalent.  
Indeed, explicit isomorphisms between various models have been demonstrated, for instance 
between the differential K-theory group of \cite{HS} and \cite{FL} in \cite{Kl}, 
between Lott's $\R/\Z$ K-theory and Lott-Freed differential K-theory in the latter \cite{FL},  
the relation between Bunke-Schick differential K-theory and Lott(-Freed) 
differential K-theory is given in \cite{Ho2}, and the 
 isomorphism between Simons-Sullivan \cite{SSu} and Freed-Lott \cite{FL} is given in 
\cite{Ho1}. 

The group structure of differential K-theory splits into odd and even degree parts, 
thus  the refinement preserves the grading. However, the odd part turns out to be 
more delicate than the even part. In particular,  while any two differential extensions of 
even K-theory are isomorphic by the uniqueness results of \cite{BS2},  odd K-theory 
requires extra data in order to obtain uniqueness. There are various concrete models
in the odd case, using smooth maps to the unitary group \cite{TWZ1}, via loop bundles 
\cite{HMSV}, and via Hilbert bundles \cite{GL}. Our results in both even and odd K-theory 
will, of course, not depend on the particular model chosen.

 Suppose $\E$ is a spectrum and $X$ is a space of the homotopy type of a CW-complex. 
Then there is a half-plane spectral sequence (AHSS) 
$$
E_2^{p, q} \cong H^p(X; \E^q(\ast))\;,
$$
converging conditionally to $\E^*(X)$. 
An immediate matter that we encounter in setting up the spectral sequence
which calculates the generalized differential cohomology of a smooth manifold 
$X$ is how to deal with filtrations. Classically, Maunder \cite{Mau} gave two approaches to 
any generalized cohomology theory. The first is by filtering over the $q$-skeletons  $X^q$ 
of the topological space $X$, and second by filtering over the Postnikov systems of spaces 
$Y_{q}$, which are the layers of an $\Omega$-spectrum associated to the cohomology 
theory. Maunder also gives an isomorphism 
between the two approaches. While we expect this to be the case in the differential 
setting, the proof might require considerable work. Hence, we leave this as an
open problem. Maunder sets  up his construction in the simplicial complex setting,
which is equivalent to setting up in the CW-complex setting as 
the geometric realization of a simplicial set is 
a CW-complex.
Simplicial and {\v C}ech spectral sequences are discussed by 
 May and Sigurdsson in \cite{MS} (Ch 22). 

We will prefer the filtration of the spaces/manifolds rather than of the corresponding 
spectra, as this will naturally bring out the geometry desired in the smooth setting. 
  We first would like to replace a topological space with skeletal filtration 
   by a smooth manifold and then view this manifold as 
 stack. Hence, in doing this, we  need an analog of a skeleton in stacks. 
This will be done via {\v C}ech resolution  of smooth spaces,
 and the replacement of skeletons of a space $X$ will be  
 the various intersections of open sets covering the smooth manifold $X$.

We will use ${\rm diff}(\Sigma^n \E, {\rm ch})$ to denote the differential refinement
in degree $n$ of a cohomology theory $\E$. This was the notation used in \cite{HS} and 
carries more data than other notation, such as 
$\E(n)$. It also avoids possible confusion with other notations, e.g. when dealing with 
Morava K-theory $K(n)$ at chromatic level $n$. 
The axiomatic approach is very useful for characterizing a smooth cohomology theory,
but one still needs the model of \cite{HS} for actually constructing examples of such 
smooth spectra. We will be using features of two main approaches at once, namely 
from \cite{HS} with $I: {\rm diff}(\Sigma^n \E, {\rm ch}) \to \underline{\E}$ and 
from \cite{BNV} \cite{Urs} with $I: \E \to {\Pi} \E$. Note that $\underline{\E}$ is not 
discrete while ${\Pi} \E$ is, but both are equivalent as smooth spectra 
$\underline{\E} \simeq {\Pi}\E$. This essentially boils down to the fact that since 
${ \Pi}\E$ is locally constant, the underlying theory satisfies ${ \Pi}\E^*(U)={ \Pi}\E^*(\ast)$ 
on contractible open sets. On the other hand, the homotopy invariance of the theory 
$\underline{\E}$ implies the same thing: namely, $\underline{\E}(U)\simeq \underline{\E}(\ast)$,
for a contractible $U$.  These relationships are discussed in further detail in \cite{BNV}.

We will be interested in how the differentials look like in our spectral sequences. 
One might a priori suspect that the differentials in the refined theories should 
at least loosely be connected to the differentials of the underlying topological theory. 
We will make this precise below, and so it seems appropriate to understand the 
form and structure of the differentials in the topological case. To illustrate the point, 
we will focus on what might perhaps be the most prominent example, namely the first differential 
$d_3: H^*(X, K^0(*)) \to H^*(X, K^0(*))$ in complex topological K-theory $K(X)$ of a topological space $X$. 
This is given by $Sq^3_\Z$ \cite{AH} \cite{AH2}.  There are exactly two stable cohomology 
operations $H^*(X; \Z) \to H^{*+3}(X; \Z)$, since $H^{n+3}(K(\Z, n))=\Z/2$ for $n$ sufficiently large. 
One of these is zero and  the other is $\beta \circ Sq^2 \circ \rho_2$, where $\beta$ 
is the Bockstein associated to the sequence $\Z \overset{\times 2}{\longrightarrow} \Z 
\overset{\rho_2}{\longrightarrow} \Z_2$ with $\rho_2$ denoting both the mod 2 reduction
and its effect on cohomology with these as coefficients, i.e. $\rho_2: H^{i}(X; \Z) \to H^{i}(X; \Z/2)$.

The above class, which is a priori in mod 2 cohomology, turned out to be a class in integral cohomology. 
One could work at any prime \cite{AH2}, by noting the following (see e.g. \cite{FFG} or \cite{Ha}).
For any class $x \in H^n(X; \Z/p)$, and with $\beta_p$ the Bockstein associated with the 
sequence $\Z_p \overset{\times p}{\to} \Z_{p^2} \overset{\rho_p}{\longrightarrow} \Z_p$,  
the elements $\beta_p(x)$ is an integral class in $H^{n+1}(X; \Z/p)$, i.e. it belongs to the 
image of the reduction homomorphism $\rho_p: H^{n+1}(X; \Z) \to H^{n+1}(X; \Z/p)$. 
This can be used to prove the integrality of the class $d \in H^3(K(\Z/p, 2); \Z/p)$ as follows 
(see \cite{FFG}). The cohomology Serre spectral sequence for the path-loop fibration 
$\Omega K(\Z, 2) \to PK(\Z, 3) \to K(\Z, 3)$ gives that $H^*(K(\Z, 3); \Z/p)$ has a single additive 
generator $\overline{d}$ in dimension $\leq 2p$. Now we have a map $\beta: K(\Z/p, 2) \to K(\Z, 3)$ 
such that 
$
\beta^*(\overline{d})=d\in H^3(K(\Z/p, 2); \Z/p)
$,
constructed via the Serre spectral sequence of the 
path-loop fibration $K(\Z/p, 1) \to PK(\Z/p, 2) \to K(\Z/p, 2)$.
The map $\beta$ induces a map of loop spaces
which are also Serre fibrations 
$$
\xymatrix{
K(\Z/p, 1) \ar[r] & PK(\Z/p, 2) \ar[d] & PK(\Z, 3)\ar[d] & K(\Z/p, 2) \ar[l]
\\
& K(\Z/p, 2) \ar[r] & K(\Z, 3)\;. &
}
$$
The induced homomorphism on the special sequences sends 
$\overline{d}$ to $d$ by the construction of $\beta$. Now we have
$H^3(K(\Z/p, 2); \Z/p)=\Z/p$ hence 
$d$ is contained is contained in the image of the homomorphism 
$\rho_p: H^3(K(\Z/p, 2); \Z) \to H^3(K(\Z/p, 2); \Z/p)$. 
Therefore $d$ is an integral class. This is attractive as it makes it readily 
amenable to differential refinement.

Such statements, and generalizations to other primes and to other generalized 
cohomology theories, can be made at the level of spectra (see e.g. \cite{Sch}). 
The first nontrivial $k$-invariant of connective complex K-theory spectrum 
$ku$ is a morphism ${\bf k}_2(ku) \in H^2(H\Z, \Z)$, which is equal to 
$\beta \circ Sq^2$, where $\beta: H\Z/2 \to \Sigma(H\Z)$ is the 
Bockstein operator associated to the extension 
$\Z \overset{\times 2}{\longrightarrow} \Z \longrightarrow \Z/2$, and 
$Sq^2_\Z$ is the pullback of the Steenrod operation 
$Sq^2 \in H^2(H\Z/2, \Z/2)$ along the projection morphism 
$\rho_2: H\Z \to H\Z/2$ given by mod 2 reduction. 
Since $ku$ is a symmetric ring spectrum then,  by \cite{Sch} (Prop. 8.8),
the $k$-invariants are derivations. The only derivations (up to units) 
in the mod $p$ Steenrod algebra ${\cal A}_p$ are the Milnor primitives
$Q_n \in H^{2p^n-1}(H\Z/p , \Z/p)$. At the lowest level we have
$Q_0=\beta_p$ the mod $p$ Bockstein, and the others
are realized as $k$-invariants of symmetric spectra, the 
connective Morava K-theory spectra $k(n)$. That is we have
$Q_n={\bf k}_{2p^n-2}(k(n))$. We will consider  refinements of integral
lifts of these. 

The classical AHSS collapses already at the first page if the generalized cohomology 
theory is rational. In fact, it can be shown that for any reasonably behaved spectrum 
like all the ones we consider, all the differentials in the AHSS are
torsion, i.e. are zero when rationalized (see \cite{Ru} Cor. 7.12).  
The differentials in the AHSS in the topological case are analyzed by 
systematically by Arlettaz \cite{Ar}. Using the structure of the integral homology of 
the Eilenberg-MacLane spectra, it is proved there that for any connected space $X$ 
there are integers $R_r$ such that $R_rd_r^{s,t}$ for all $r \geq 2$, $s$, $t$. 
Some aspects of this general feature will continue to hold in the differential 
setting. Form a homotopy point of view there is not much difference between the 
localizations at $\R$ and at $\QQ$. However, from a geometric point of view there is 
a considerable difference. Nevertheless, we will still use the term ``rationalize" when 
we discuss localization at $\R$, as customary in the homotopy theory literature. 
We stress that the distinction is needed in certain geometric settings (see \cite{GM}), 
but it will not be an issue for us in this paper.

Note that although the differential cohomology diamond, i.e. the diagram that characterizes such theories 
(see Remark \ref{diff diamond}), certainly detects torsion classes in the flat part of the theory, it does not distinguish between torsion at various primes. As a by-product, our analysis can be seen as a systematic method for addressing $p$-primary torsion in differential theories. In \cite{GS2} we found that the 
Deligne-Beilinson squaring operation admits lower degree operations refining the Steenrod squares. 
We have the familiar pattern
$$
{\rm DD}, \widehat{Sq}^1,  \widehat{Sq}^2,  \widehat{Sq}^3,  \cdots {\rm DD}^2, \cdots \;,
$$
where ${\rm DD}$ is the Dixmier-Douady class: a non-torsion differential cohomology operation. The refined 
squares $\widehat{Sq}^{2k+1}$, as the classical squares $Sq^{2k}$, are operations that are 2-torsion. 
In this paper we get $\widehat{Sq}^{2k+1}$ as we expect, but also differentials 
$d_{2m}$ at lowest degree for every $m$,
\(
d_{2m}: \prod_{k}\Omega^{2k}(M) \longrightarrow H^{2m}(M; U(1))\;.
\label{eq d2m}
\) 
We consider this as a cohomology operation, which can be viewed as first projecting on to the homogeneous component ${\rm ch}_{2m}$ of the Chern-character. A $U(1)$-valued {\v C}ech-cocycle is obtained by restricting to $2m$-fold intersections of an open cover, pairing with an appropriate simplex of degree $2m$ and exponentiating (This will be spelled out in detail in section \ref{Sec App}). If indeed the form ${\rm ch}_{2m}$ arises as the curvature of a bundle, it must represent a closed form with \emph{integral} periods. The differential $d_{2m}$ can therefore be understood as the obstruction to this condition. 
Similar results hold for the odd part, i.e. for differentially refined $K^1$-theory, where 
the refined Steenrod square takes the same form as in differential $K^0$-theory , while the 
differentials arising from forms -- the analogues of those in \eqref{eq d2m} -- 
are now of odd degrees. 

The paper is organized as follows. In Sec. \ref{Sec Smooth} we start by carefully setting up the background 
in smooth and differential cohomology, preparing the scene for our constructions. In particular, in 
Sec. \ref{Sec stack} we adapt abstract general results on stacks (or simplicial sheaves) 
to our context and spell out specific definitions and constructions that will be useful for
us in later sections; more general and comprehensive accounts can be found in \cite{Ja} \cite{Lu} \cite{Urs}. 
Then in Sec. \ref{Sec Diff} we take the approach to differential cohomology
that allows for a direct generalization. Our main constructions will be in Sec. \ref{Sec AHSS}, 
and in particular in Sec. \ref{Sec Cech} we provide the filtration via {\v C}ech resolutions
and then construct the AHSS for smooth spectra in Sec. \ref{Sec Morph} and compare to
the AHSS of the underlying topological theory. This refinement 
will depend on whether the degree is positive, negative, or zero. Then we explore the compatibility 
of the differentials with the product structure in Sec. \ref{Sec Product}. 

Having given the main construction, our main applications of the general spectral sequence   
to various differential  cohomology theories will be presented in Sec. \ref{Sec App}.  
  The construction is general enough to apply to any structured cohomology theory
whose coefficients are known. We will explicitly emphasize three main examples: ordinary differential 
cohomology, differential K-theory,  and a differential version of integral Morava K-theory that we introduce. 
As a test of our method, in Sec. \ref{Sec Coh} we recover the usual hypercohomology spectral sequence for the 
Deligne complex (see \cite{Bry}, \cite{EV} Appendix), and we do so for  manifolds, then products of these, 
and then more generally for smooth fiber bundles. Then the AHSS for K-theory is generalized in 
Sec. \ref{Sec K} to differential K-theory, where the differential involve refinements of Steenrod squares, 
in the sense of \cite{GS2}, as well as operations on forms, as indicated above around expression
 \eqref{eq d2m}. 
We also show that the odd case, i.e. smooth extension of $K^1$, leads to a similar construction, 
but with the differentials now involving odd forms. 
Then in Sec. \ref{Sec Kn} we first introduce a refinement of the integral 
form of Morava K-theory, discussed in \cite{KS1} \cite{S1} \cite{SW}, and then 
characterize the corresponding differentials, which turn out to have a similar pattern as in 
K-theory, where the operation that gets refined is the Milnor primitive $Q_n$ encountered above.
We end with an application to an example from M-theory and string theory.

\paragraph{Notation.} We have the following morphism that we will use repeatedly throughout.
Denote by  $\rho_p:\ZZ\to \ZZ/p$  the mod $p$ reduction on coefficients with corresponding 
morphism with the same notation on the cohomology groups with these as coefficients. 
 We will denote by $\beta$, $\beta_p$, and $\tilde{\beta}$ the 
Bockstein homomorphisms associated with  the 
coefficient sequences
$$0\to \ZZ \to \RR \overset{\exp}{\longrightarrow} U(1)\to 0\;,$$
$$0\to\ZZ/p \overset{\times p}{\longrightarrow}\ZZ/{p^2} \overset{\rho_p}{\longrightarrow} \ZZ/p \to  0\;,$$
$$0\to\ZZ \overset{\times p}{\longrightarrow}\ZZ \overset{\rho_p}{\longrightarrow} \ZZ/p \to  0\;,$$
respectively. 
We will let $\Gamma_2:\ZZ/2\into U(1)$ denote the representation as the square roots of unity, also 
with  $\Gamma_2:H^n(-;\ZZ/2)\to H^n(-; U(1))$ the induced map on cohomology. 
 We will also use more refined Bockstein homomorphisms associated with spectra, 
 and these will be defined as we need them.

\section{Smooth cohomology} 
\label{Sec Smooth}

\subsection{Smooth cohomology and stable category of smooth stacks}
\label{Sec stack}

In this section we adapt abstract general results on stacks (or simplicial sheaves) to our context 
and spell out specific definitions and constructions that will be useful for us in later sections. The 
interested reader can find more general and comprehensive accounts in \cite{Ja} \cite{Lu} \cite{Urs}. 
For the reader who is more interested in the applications to differential cohomology theories, this section 
can be skipped. However, we would like to emphasize that although the language used in this section is 
rather abstract, the generality gained from this formalism is far reaching and allows this machinery to be 
used for a wide variety theories, beyond just differential cohomology theories. 

Essentially, the axioms characterizing a smooth cohomology theory are not much different from the axioms 
characterizing usual cohomology theories. The big difference is where the theory takes place. More precisely, 
we want to consider homotopical functors on the category of pointed \emph{smooth stacks} 
$Sh_{\infty}(\cartsp)_+$ with $\cartsp$ the category of Cartesian spaces, rather than the category 
of pointed topological spaces $\Top_+$. Let $\ab_{gr}$ be the category of graded abelian groups.

\begin{definition} (Smooth cohomology) Let $\E^*:Sh_{\infty}(\cartsp)_+^{\rm op} \to \ab_{gr}$ be a functor 
satisfying the following axioms:
\begin{enumerate} 
\item {\rm (Invariance)} $\E^*$ sends equivalences to isomorphisms.
\item {\rm (Additivity)} For small coproducts (i.e. ones forming sets) of pointed stacks, 
$\bigvee_{\alpha}{X}_{\alpha}$, we have 
$$
\E^*\Big(\bigvee_{\alpha} {X}_{\alpha}\Big)=\prod_{\alpha}\E^*({X}_{\alpha})\;.
$$
\item {\rm (Mayer-Vietoris)} For any homotopy pushout of pointed stacks,
$$\xymatrix{
Z\ar[r]\ar[d] & {Y}\ar[d]
\\
{X}\ar[r] & {X\cup_{Z}Y}\;,
}
$$
the induced sequence 
$$\E^*({X\cup_{Z}Y})\to \E^*({X})\oplus \E^*({ Y})\to \E^*({Z})$$
is exact.
\item {\rm (Suspension)} For any stack ${X}$, there is an isomorphism $\E^{n+1}(\Sigma {X})\simeq \E^n({X})$.
\end{enumerate}
Then we call $\E^*$ a {\rm smooth cohomology theory}.
\end{definition}

\begin{remark}
Note that the Mayer-Vietoris axiom implies the usual Mayer-Vietoris sequence. Indeed, let $M$ be a manifold 
and let $V$ be a local chart of $M$. Let $U$ be an open set such that $\{U,V\}$ is a cover of $M$. Then 
the strict pushout
$$\xymatrix{
U\cap V\ar[r]\ar[d] & V\ar[d]
\\
U\ar[r] & U\cup V
}
$$
is actually a homotopy pushout. We can equivalently write this diagram as a homotopy coequalizer
$$
\xymatrix{
U\cap V \ar@<.1cm>[r] \ar@<-.1cm>[r] & U\coprod V \ar[r] & U\cup V\;,
}
$$
in which the homotopy cofiber of the second map can be identified with $\Sigma U\cap V$. By iterating this argument and applying $\E^*$ to the the resulting diagram one obtains the long exact sequence
$$
\hdots \to \E^{*}(U\cap V) \to \E^*(M)\to \E^*(U)\oplus \E^*(V)\to \E^{*+1}(U\cap V) \to \hdots \;,
$$
which is the familiar Mayer-Vietoris sequence.
\end{remark}
The above axioms can be taken as a generalization of the Eilenberg-Steenrod axioms (see \cite{Ad} \cite{Hi}), 
where the Mayer-Vietoris axiom subsumes both the excision axiom and the long-exact sequence axiom. 
 It is interesting to note that the axioms \emph{do not} require homotopy invariance. Namely, if two 
 manifolds $M$ and $N$ are \emph{homotopic}, they may fail to be equivalent as stacks. In fact, 
 an equivalence of stacks requires, in particular, that for every sheaf $F$ (embedded as a stack), 
 we have an isomorphism
$$
F(N)\simeq \pi_0\map(N, F)\simeq \pi_0\map(M, F)\simeq F(M)\;.
$$
In particular, we can take the sheaf of smooth $\RR$-valued functions on a 
manifold. Then if every homotopy equivalence 
$f:M\to N$ induced an equivalence of stacks, we would have an induced isomorphism
$$
f^*:C^{\infty}(N;\RR) \to C^{\infty}(M;\RR)\;.
$$
Taking $N=\ast$ and $M=\RR^n$ immediately gives a contradiction. On the other hand, every equivalence 
of stacks does produce a weak homotopy equivalence of geometric realizations. To see this, simply note that 
the geometric realization functor
$$
\Pi:Sh_{\infty}(\cartsp) \to \sset\;,
$$
 being a Quillen functor, has a derived functor by Ken Brown's Lemma \cite{Br}. It therefore 
 preserves weak equivalences between fibrant objects. But these objects are exactly those 
 that satisfy descent, namely stacks, (e.g. manifolds) \cite{Urs} \cite{Dug}.

\begin{remark}
Given a smooth cohomology theory $\E^*$, we always get a presheaf of graded abelian groups on 
the site $\cartsp$ by precomposing with the Yoneda embedding:
$$
\xymatrix{
\E^*:\cartsp~ \ar@<-.2ex>@{^{(}->}[r]^{Y} & Sh(\cartsp)~\ar@<-.2ex>@{^{(}->}[r]^{{\rm sk}_0} &
 Sh_{\infty}(\cartsp)\ar[r]^-{+} &Sh_{\infty}(\cartsp)_+ \ar[r]^-{\E^*}& \ab_{\rm gr}
 }\;,
$$
where ${\rm sk}_0$ embeds a sheaf as a discrete simplicial sheaf.
We will use this fact later in the construction of the spectral sequence in theorem \ref{main theorem}.
\end{remark}

Just as all cohomology theories are representable by $\Omega$-spectra, via Brown representability, all smooth 
cohomology theories are representable by smooth spectra. This follows from the version of Brown representability
 formulated by Jardine in \cite{Ja} applied to the stable homotopy category of smooth stacks. We will quickly review
  the basic properties of this category (see \cite{Lu} \cite{Ja2}) to establish where our objects of interest live. 

We first recall some operations on stacks that are counterparts to standard operations on 
topological spaces. Let $X$ and $Y$ be two pointed stacks. 
\begin{enumerate}[label=(\roman*)]
\item The wedge product $X \vee Y$ is defined via the pushout diagram 
$$
\xymatrix{
Y \ar[r] & Y \vee X \\
\ast \ar[u] \ar[r] & X\;. \ar[u]
}
$$
\item The smash product $ X \wedge Y$ is defined as 
the quotient 
$X \wedge Y:= X \times Y/X \vee Y$ of the Cartesian product by the wedge product. 

\item The suspension $ \Sigma X$ is defined via the homotopy pushout diagram
$$
\xymatrix{
X \ar[d] \ar[r] & \ast \ar[d] \\
\ast \ar[r] & \Sigma X\;.
}
$$
\item The looping, i.e. loop space, $\Omega X$ is defined via the homotopy pullback
$$
\xymatrix{
\Omega X \ar[d] \ar[r] & \ast \ar[d] \\
\ast \ar[r] & X\;.
}
$$
\end{enumerate}

\begin{definition} 
We define the {\rm stabilization ${\rm Stab}(Sh_{\infty}(\cartsp)_+)$ of smooth pointed stacks} to be the following
 category:
\begin{list}{$\circ$}{}  
\item The objects of ${\rm Stab}(Sh_{\infty}(\cartsp)_+)$ are sequences of pointed stacks 
$$\{ \E_n\}\subset Sh_{\infty}(\cartsp)_+,\ \ n\in \ZZ$$
equipped with maps
$\sigma_n:\Sigma \E_n\to {\E}_{n+1}$.
\item The morphisms between $\E$ and ${\cal F}$ are 
defined to to be the  levelwise morphisms 
${\E}_n \to {\cal F}_n$,
commuting with the $\sigma_n$'s.
\end{list}
\end{definition}
This category carries a stable model structure given by first taking the projective model structure on sequences of 
stacks and then performing Bousfield localization with respect to stable weak equivalences in the usual way. This process is described in detail in \cite{Ja} \cite{Lu} \cite{Ja2} and we summarize the relevant results found there.
 The category ${\rm Stab}(Sh_{\infty}(\cartsp)_+)$ admits a stable, closed, simplicial model structure in which
\begin{list}{$\circ$}{}  
\item The weak equivalences are \emph{stable} weak equivalences. That is, a morphism of smooth spectra 
$f:{ \E}_{\bullet}\to {\cal F}_{\bullet}$ is a weak equivalence if and only if it induces a weak equivalence
$$
Q(f):\lim_{i\to \infty} \Omega^i{ \E}_{n+i} \to \lim_{j\to \infty}\Omega^j{\cal F}_{n+j}\;.
$$
\item The fibrant objects are precisely the smooth $\Omega$-Spectra, that is, the sequence of stacks ${X}_{\bullet}$ whose structure maps
$$
\sigma_n:\Sigma {\E}_n\to {\E}_{n+1}
$$
induce equivalences
${\E}_n\overset{\sim}{\to} \Omega {\E}_{n+1}$.
\end{list}

\begin{remark} 
We will refer to the stable model category ${\rm Stab}(Sh_{\infty}(\cartsp)_+)$ as the category of 
\emph{smooth spectra} and denote it by
$$Sh_{\infty}(\cartsp;\mathscr{S}{\rm p}):={\rm Stab}(Sh_{\infty}(\cartsp)_+)\;.$$
\end{remark}

\begin{example}
Let $M\in Sh_{\infty}(\cartsp)_+$ be a manifold, viewed a stack and equipped with a basepoint. We 
can define the smooth spectrum $\Sigma^{\infty}M$ in the usual way, as the sequence of 
suspensions of the manifold $M$. Given a smooth 
 $\Omega$-spectrum ${\E}$, we can define a smooth cohomology theory $\E^*$, by setting 
$$\E^q(M)\simeq \pi_0\map(\Sigma^{-q}\Sigma^{\infty}M,{ \E})\;.$$
\end{example}

Differential cohomology theories are examples of the theories introduced above, although it may not be immediately 
apparent where the differential cohomology ``diamond" diagram \cite{SSu} fits into this context. In fact, it was observed 
by Bunke, Nikolaus and V\"olkl in \cite{BNV}, that the diamond provides a further characterization of \emph{all} smooth cohomology theories in terms of refinement of topological theories. This characterization happens in addition to the 
Brown representability described above, and happens only when the category of stacks exhibits so-called 
\emph{cohesion}. We now review the properties of the cohesive structure on smooth stacks \cite{Urs} that we need, 
along with the characterization of smooth cohomology theories described in \cite{BNV}. 
It is shown in \cite{Urs} that the category $Sh_{\infty}(\cartsp)$ admits a quadruple $\infty$-categorical adjunction 
$(\Pi \dashv {\rm disc} \dashv \Gamma \dashv {\rm codisc})$ 
\(
\xymatrix{
Sh_{\infty}(\cartsp)\ar@<-.5em>[rr]^<<<<<<<<<<{\Gamma} \ar@<1.5em>[rr]^<<<<<<<<<<{\Pi} && \sset \ar@<1.5em>[ll]_<<<<<<<<<<{\rm codisc} \ar@<-.5em>[ll]_<<<<<<<<<<{\rm disc} 
}
\label{cohesion}
\)
where $\Pi$ preserves finite $\infty$-limits and the functors ${\rm disc}$ and ${\rm codisc}$ are fully faithful.  

One implication of this is that $\sset$ embeds into $Sh_{\infty}(\cartsp)$ as an $\infty$-subcategory
 in two different ways, one reflective, the other reflective and
coreflective. From the reflectors one can produce two monads and one comonad defined as follows:
$$
{ \Pi}:=\Pi\circ {\rm disc},\qquad
 \flat:={\rm disc}\circ \Gamma, \qquad 
  \sharp:={\rm codisc}\circ \Gamma\;.
  $$
These monads fit into a triple $\infty$-adjunction $({ \Pi}\dashv \flat \dashv \sharp)$ which is called a \emph{cohesive} 
adjunction. 

\begin{remark}
Each monad in the cohesive adjunction picks out a different part of the nature of a smooth stack. This nature is 
perhaps best exemplified by how the adjoints behave on smooth manifolds (viewed as stacks). More precisely, 
if $M$ is a smooth manifold
then, for instance, 
\begin{enumerate}[label=(\roman*)]
\item 
the comonad $\flat$ takes the underlying set of points of the manifold and then embedds this set back into 
stacks as a discrete object. This functor therefore misses the smooth structure of the manifold and treats it 
instead as a discrete object. 

\item The monad ${ \Pi}$ essentially takes the singular nerve of the maniflold using \emph{smooth} paths 
and higher smooth simplices on the manifold. It therefore retains the geometry of the manifold and ``knows" 
that the points of the manifold ought to be connected together in a smooth way. 
\end{enumerate}
\end{remark} 

\noindent The following observation on lifting from simplicial sets to spectra 
is known (\cite{Urs}, Prop. 4.1.9), but we supply a proof for completeness. 

\begin{proposition} 
\label{Propquad}
The $\infty$-adjunction \eqref{cohesion} 
lifts to an $\infty$-adjunction
$$
\xymatrix{
Sh_{\infty}(\cartsp;\mathscr{S}{\rm p})\ar@<-.5em>[rr]^<<<<<<<<<<{\Gamma^s} \ar@<1.5em>[rr]^<<<<<<<<<<{\Pi^s} && \mathscr{S}{\mathrm p} \ar@<1.5em>[ll]_<<<<<<<<<<{{\rm codisc}^s} \ar@<-.5em>[ll]_<<<<<<<<<<{{\rm disc}^s} 
}$$
on the stable $\infty$-category of smooth spectra. Moreover, the adjoints satisfy the same condition as the 
$\infty$-adjunction \eqref{cohesion} does.
\end{proposition}
\theproof
The category of smooth stacks is presented by the combinatorial simplicial model category 
$$
Sh_{\infty}(\cartsp)=[\cartsp,\sset]_{{\rm loc},{\rm proj}}\;,
$$
where ${\rm loc}$ denotes the Bousfield localized model structure at the maps out of {\v C}ech nerves. The quadruple 
adjunction is presented by Quillen adjoints $(\Pi\dashv {\rm disc} \dashv \Gamma \dashv {\rm codisc})$ \cite{Urs}. We
 need to show that this adjunction holds on the stable model category of smooth spectra. The adjunction immediately 
 gives an underlying categorical adjunction by simply applying the functors degree-wise. In the projective model structure, the right adjoints are Quillen by definition and the closed model axioms imply that the left adjoints are also Quillen. 

Now the functors (in the global model structure on $\mathscr{S}{\rm p}$) ${\rm disc}$ and ${\rm codisc}$ both 
preserve homotopy limits. Hence for a local weak equivalence $f:\E \to {\cal F}$ of spectra, we have
\bea
\lim_{i\to \infty} \Omega^i{\rm disc}(\E)_{n+i}  &\simeq & {\rm disc}\big(\lim_{i\to \infty} \Omega^i {\cal F}_{n+i}\big) 
\\
& \simeq & {\rm disc}\big(\lim_{j\to \infty}\Omega^j {\cal F}_{n+j}\big)
\\
&\simeq & \lim_{j\to \infty}\Omega^j{\rm disc}({\cal F})_{n+j}
\eea
and ${\rm disc}(f)$ induces a weak equivalence $Q({\rm disc}(f))$. Hence, ${\rm disc}(f)$ is a weak equivalence. In the
 same way, ${\rm codisc}$ preserves local weak equivalences. It follows by the basic properties of Bousfield localization
  that ${\rm disc}$ and ${\rm codisc}$ are right Quillen adjoints. Again, by the axioms of a closed model category, it follows that the entire adjunction holds as Quillen adjunction of stable model categories.
\endofproof
\begin{remark} The proof of the previous proposition implies that both ${\rm disc}$ and ${\rm codisc}$ preserve 
$\Omega$-spectra. However, ${ \Pi}$ and $\Gamma$ need not take $\Omega$-spectra to $\Omega$-spectra. 
This problem can be remedied by taking $\Pi^s$ (or $\Gamma^s$) to be the composite of $R\circ \Pi$ 
(or $R\circ \Gamma$), where $R$ is the fibrant replacement in spectra. Since $R$ defines a left $\infty$-adjoint
 to the inclusion of fibrant objects (and preserves finite $\infty$-limits), we will still have an adjunction at the level of
  $\infty$-categories (although this is not presented by a Quillen adjunction). 
\end{remark}

As in the case of smooth stacks, the quadruple adjunction in the Proposition \ref{Propquad} produces adjoint monads 
$({\Pi}^s\dashv \flat^s \dashv \sharp^s)$ exhibiting \emph{stable} cohesion. The main observation in \cite{BNV}, 
recast in the cohesive setting in \cite{Urs}, is the following. Let $j:\flat^s\to {\rm id}$ be the counit of the 
comonad $\flat^s$, and let $I:{\rm id}\to { \Pi}^s$ be the unit of the monad ${ \Pi}^s$. Let 
${\E}\in Sh_{\infty}(\cartsp;\mathscr{S}{\rm p})$ be a smooth spectrum. Then ${\E}$ sits inside a hexagon diagram
\(\hspace{-.8cm}\label{spectra diamond}
\xymatrix@C=.5pt @!C{
&{\rm fib}(\eta)({ \E}) \ar[rd]\ar[rr] & & {\rm cofib}(\epsilon)({ \E}) \ar[rd] &
\\
\Sigma^{-1}{ \Pi}^s{\rm cofib}(\epsilon)({ \E})\ar[ru]\ar[rd] & & {\E}\ar[rd]^-{I} \ar[ru]& &  { \Pi}^s{\rm cofib}(\epsilon)({\E})\;,
\\
&\flat^s {\E}\ar[ru]^{j}\ar[rr] & & { \Pi}^s{ \E}\ar[ru] &
}
\)

\vspace{2mm}
\noindent where the diagonals are fiber sequences (by definition), the top and bottom sequences are fiber 
sequences, and the two squares in the hexagon are homotopy Cartesian, i.e.
both are homotopy pullback squares and hence homotopy pushout (via the equivalence of the
two in the stable setting). The latter property is key, because it is a homotopy Cartesian square,
as on the right of the hexagon, which Hopkins-Singer \cite{HS} took as the definition of
differential cohomology (for a specific choice of the object of
differential forms). Bunke-Nikolaus-V\"olkl \cite{BNV} observed that by the
hexagon, {\it every} smooth spectrum satisfies this kind of Hopkins-Singer
definition, if one just allows more general objects of differential
forms, which is the object ${\rm cofib}(\epsilon)(\E)$ in our notation above.

It often happens in practice that the smooth spectra ${\rm fib}(\eta)({\E})$ and ${\rm cofib}(\epsilon)({ \E})$ contain 
no information away from degree 0. In particular, it often happens that for $n>0$,
\begin{eqnarray}
\pi_{n} {\rm Map}\big(M, {\rm cofib}(\epsilon)({ \E}) \big) &\simeq & 0\;,
\label{n pm1}
\\
\pi_{-n} {\rm Map}\big(M, {\rm fib}(\eta)({ \E})\big) &\simeq & 0\;.
\label{n pm2}
\end{eqnarray}
In this case the ${ \E}$-cohomology of a manifold can be calculated as either the flat cohomology
or the underlying topological cohomology in all degrees but 0. This is summarized as the following result.

\begin{proposition}\label{diff-cohomology-cal}
 Let ${\E}$ be a smooth spectrum such that \eqref{n pm1} and \eqref{n pm2} are satisfied.
Then the ${\E}$-theory of a manifold $M$ is given by
$${ \E}^n(M):=\left\{\begin{array}{ccc}
({ \Pi}^s \E)^n(M) && n>0,
\\
\\
(\flat^s \E)^n(M) && n<0, 
\end{array}\right.$$
(where in degree 0, $\E(M)$ is already characterized by the diamond \eqref{spectra diamond}).
\end{proposition}
\theproof
Since the diagonals of the diamond are fiber sequences, they induce long exact sequences in cohomology. 
Let $n$ be a positive integer. The sequence
$$\flat^s {\E} \to { \E} \to {\rm cofib}(\epsilon)({ \E})$$
gives the section of the long sequence
$$
\pi_{n+1}\map(M,{\rm cofib}(\epsilon)({\E})) \to \flat^s { \E}^{-n}(M) \to { \E}^{-n}(M) \to \pi_{n}\map(M,{\rm cofib}(\epsilon)({ \E}))\;.
$$
By assumption the leftmost and rightmost groups are $0$. We therefore have an isomorphism
$$
(\flat^s { \E})^{-n}(M)\simeq  { \E}^{-n}(M)\;.
$$ 
Similarly, the sequence
$$
{\rm fib}(\eta)({\E}) \to { \E} \to {  \Pi}^s \E\;,
$$
gives the long sequence
$$
\pi_{-n}\map(M{\rm fib}(\eta)({ \E})) \to { \E}^n(M) \to ({ \Pi}^s \E)^n(M) \to \pi_{-n-1}\map(M,{\rm fib}(\eta)({ \E}))\;,
$$
and again we get the desired isomorphism.
\endofproof

\subsection{Differential cohomology and differential function spectra}
\label{Sec Diff}

The main applications we have in mind, as we indicated in the Introduction, 
 concern \emph{differential cohomology theories}. In this section we review some of the concepts established in 
 \cite{Bun} \cite{BNV} \cite{Urs} (which generalize \cite{SSu}), adapted to our context. 

\begin{definition}
Let $\E^*$ be a cohomology theory. A {\rm differential refinement} $\widehat{\E}^*$ of $\E^*$ consists of the 
following data:
\begin{enumerate} 
\item A functor $\widehat{\E}^*:Sh_{\infty}(\cartsp_+)^{\rm op} \to \mathscr{A}\mathrm{b}_{gr}$;
\item Three natural transformations:
\begin{enumerate}
\item {\rm Integration:} $I:\widehat{\E}^*\to \E^*$;
\item {\rm Curvature:} $R:\widehat{\E}^*\to Z_*\left(\Omega^*\otimes \E^*(\ast)\right)$;
\item {\rm Secondary Chern character:} $a:\Omega^*\otimes\E^*(\ast)[1]/{\rm im}(d)\to \widehat{\E}^*$;
\end{enumerate}
\end{enumerate}
such that the following axioms hold:
\begin{list}{$\circ$}{} 
\item {\rm (Chern-Weil).} We have a commutative diagram 
$$
\xymatrix{
\widehat{\E}^*\ar[r]^-{R}\ar[d]^{I} & Z_*\left(\Omega^*\otimes\E^*(\ast)\right)\ar[d]^q
\\
\E^*\ar[r]^-{{\rm ch}} & H_*\left(\Omega^*\otimes\E^*(\ast)\right)\;,
}
$$
where ch is the Chern character map. 
\item {\rm (Secondary Chern-Weil).} We have a commutative diagram
$$
\xymatrix{
\Omega^*\otimes\E^*(\ast)[1]/{\rm im}(d) \ar[rr]^-{d} \ar[dr]_-{a} &&  Z_*\left(\Omega^*\otimes\E^*(\ast)\right)
\\
& \widehat{\E}^*\ar[ur]_-{R} & 
}
$$
and an exact sequence
$$
\hdots \to \E^*[1] \to \Omega^*\otimes\E^*(\ast)[1]/{\rm im}(d)\to \widehat{\E}^* \to \E^* \to \hdots\;.
$$
\end{list}
\end{definition}

Note that in item {\it Chern-Weil} above,  $H_*(\Omega^* \otimes \E^*(*))$ appears as 
 the codomain of the Chern character. As explained in \cite{BNV}, this becomes a
locally constant stack equivalent to just the locally constant stack on
the rationalization of $\E^*$, i.e., ch is equivalent to
${\rm ch} :  \E^* \to \E^* \wedge H\R$ (or $M\R$).

\begin{remark}\label{diff diamond}
The above characterization can ultimately be summarized by saying that differential cohomology fits into an 
exact diamond 
$$\hspace{-1cm}
\xymatrix @C=5pt @!C{
 &\Omega^*\otimes\E^*(\ast)[1]/{\rm im}(d) \ar[rd]^{a}\ar[rr]^{d} & &  Z_*\left(\Omega^*\otimes\E^*(\ast)\right)\ar[rd] &  
\\
\E^{*-1}\otimes \RR\ar[ru]\ar[rd] & & {\widehat{\E}^*}\ar[rd]^{I} \ar[ru]^{R}& &  \E^{*}\otimes \RR \;,
\\
&\E_{\RR/\ZZ}^{*-1} \ar[ru]\ar[rr]^{\beta_{\E}} & & \E^* \ar[ru]^{\rm ch} &
}
$$
where the diagonal, top and bottom sequences are all part of long exact sequences. The bottom sequence 
is obtained by observing that the cofiber of the rationalization map is an $MU(1)$ (Eilenberg-Moore spectrum),
where we identify $\RR/\ZZ$ with $U(1)$ throughout. 
That is, we have a cofiber sequence involving the unit map from the sphere spectrum $\mathbb{S}=M\Z$
$$
\mathbb{S}\to M\RR\to MU(1)\;.
$$
Smashing on the left with the theory $\E$, we obtain a ``Bockstein sequence"
$$\E \to \E\wedge M\RR \to \E\wedge MU(1)\overset{\beta_{\E}}{\longrightarrow} \Sigma\E\;.$$
We define the flat theory as
$$
\E_{U(1)}:=\E\wedge MU(1)
$$
and the rational theory as
$$
\E_{\RR}:=\E\wedge M\RR\;.
$$
\end{remark}

\begin{remark} Differential cohomology theories are a special case of smooth cohomology theories, while 
differential function spectra are a special case of smooth spectra. Thus, this section can be viewed as 
describing a special case of the previous section.
\end{remark}

Since differential cohomology theories will arise as certain homotopy pullbacks
(in Def. \ref{differential function} below), we will first need to establish the components of the pullback. 
We begin with the following lemma that can be found in \cite{Bun} (Lemma 6.10), which explains how we can 
transition from a topological cohomology theory to a smooth one, in a process 
whose direction is opposite to that of the map $I$.

\begin{lemma}
Let $\E$ be a spectrum and define the smooth presheaf of spectra $\underline{\E}$ via the assignment
\bea
Objects:&& U\mapsto \map(\Sigma^{\infty}U, \E)\;,
\nonumber\\
Morphisms:&& (f:U\to V) \mapsto (f^*:\map(\Sigma^{\infty}V, \E) \to \map(\Sigma^{\infty}U, \E))\;.
\eea
Then $\underline{\E}$ satisfies descent.
\end{lemma}
\theproof
Let $C^{\bullet}(\{U_{\alpha}\})$ denote the {\v C}ech nerve of a good open cover $\{U_{\alpha}\}$ of some 
manifold $M$. The Yoneda Lemma and basic properties of the mapping space functor imply that we have 
the following sequence of equivalences
\begin{alignat}{2}
\underline{\E}(M) &:= \map(\Sigma^{\infty}M, \E) 
\nonumber\\
&\simeq  \map(\Sigma^{\infty}{\rm hocolim}_{\Delta^{op}}C^{\bullet}(\{U_\alpha\}), \E)
\nonumber\\
&\simeq   \map({\rm hocolim}_{\Delta^{op}}\Sigma^{\infty}C^{\bullet}(\{U_\alpha\}), \E)
\nonumber\\
&\simeq {\rm holim}_{\Delta^{op}}\map (\Sigma^{\infty}C^{\bullet}(\{U_\alpha\}),\E)
\nonumber\\
&\simeq  {\rm holim}\left\{
\xymatrix{
\hdots  \ar@<-.5pc>[r] \ar@<.5pc>[r] \ar[r]  &  \ar@<-.75pc>[l] \ar@<-.25pc>[l] \ar@<.25pc>[l] \ar@<.75pc>[l]   \prod_{\alpha\beta\gamma}\map(\Sigma^{\infty}U_{\alpha\beta\gamma},\E)  \ar@<-.25pc>[r] \ar@<.25pc>[r] &  \ar@<-.5pc>[l] \ar@<.5pc>[l] \ar[l]\prod_{\alpha\beta}\map(\Sigma^{\infty}U_{\alpha\beta},\E) \ar[r] & \ar@<-.25pc>[l] \ar@<.25pc>[l] \prod_{\alpha}\map(\Sigma^{\infty}U_{\alpha}, \E)
}\right\}
\nonumber\\
&\simeq  {\rm holim}\left\{
\xymatrix{
\hdots  \ar@<-.5pc>[r] \ar@<.5pc>[r] \ar[r]  &  \ar@<-.75pc>[l] \ar@<-.25pc>[l] \ar@<.25pc>[l] \ar@<.75pc>[l]   \prod_{\alpha\beta\gamma}\underline{\E}(U_{\alpha\beta\gamma})  \ar@<-.25pc>[r] \ar@<.25pc>[r] &  \ar@<-.5pc>[l] \ar@<.5pc>[l] \ar[l]\prod_{\alpha\beta}\underline{\E}(U_{\alpha\beta}) \ar[r] & \ar@<-.25pc>[l] \ar@<.25pc>[l] \prod_{\alpha}\underline{\E}(U_{\alpha})
}\right\}\;,
\nonumber
\end{alignat}
and so $\underline{\E}$ satisfies descent.
\endofproof

The other components of the pullback we want to establish are presented by sheaves of chain complexes. 
There is a general functorial construction by which one can turn an unbounded chain complex into a spectrum, 
which we now describe (See \cite{Sh} for details). This functor is called the \emph{Eilenberg-MacLane} functor
\(
H:\mathscr{C}{\rm h} \to \mathscr{S}{\rm p}
\label{EM fun}
\)
and acts on objects as follows. 
Let $C_{\bullet}$ be an unbounded chain complex, and let $Z_n$ denote the subgroup of cycles in degree $n$. 
The functor $H$ takes $C_{\bullet}$ and forms the sequence of truncated bounded chain complexes
$$C_{\bullet}(\bullet)=\left\{
\begin{array}{l}
\left(\hdots \to C_n\to C_{n-1}\to \hdots C_1\to Z_0\right)=C_{\bullet}(0)
\\
\\
\left(\hdots \to C_n\to C_{n-1}\to \hdots C_0\to Z_{-1}\right)=C_{\bullet}(1)
\\
\\
\left(\hdots \to C_n\to C_{n-1}\to \hdots C_{-1}\to Z_{-2}\right)=C_{\bullet}(2)
\\
\vdots
\\
\left(\hdots \to C_n\to C_{n-1}\to \hdots C_{-k}\to Z_{-k}\right)=C_{\bullet}(k) 
\\
\vdots
\end{array}
\right.$$
The reason for the group of cycles appearing in degree $0$ comes from using the \emph{right} adjoint to the 
inclusion $i:\chp\to \mathscr{C}{\rm h}$ (as opposed to the left). The left adjoint simply truncates the complex 
in degree $0$, while the right adjoint truncates and then takes only the cycles in degree $0$.

Continuing with our discussion, at each level in the sequence, $H$ applies the Dold-Kan functor 
$DK: \chp \to \sset$
to the bounded chain complex in that degree. This gives a sequence of spaces
$$DK(C_{\bullet}(\bullet))=\left\{
\begin{array}{c}
DK(C_{\bullet}(0))
\\
\\
DK(C_{\bullet}(1))
\\
\\
DK(C_{\bullet}(2))
\\
\vdots
\\
DK(C_{\bullet}(k)) 
\\
\vdots
\end{array}
\right.
$$
Since $DK$ preserves looping (being a right Quillen adjoint) and equivalences (being a Quillen equivalence of 
model categories), we get induced equivalences
$$
\sigma_k:DK(C_{\bullet}(k)) \to \Omega DK(C_{\bullet}(k-1))\;,
$$
which turns $DK(C_{\bullet}(\bullet))$ into a spectrum. 
\begin{example}
Consider the unbounded chain complex $\ZZ[0]$, with $\ZZ$ concentrated in degree $0$. Then 
$$
H(\ZZ[0])\simeq H\ZZ
$$
where the right hand side denotes the Eilenberg-MacLane spectrum.
\end{example}
\begin{example}
Fix a manifold $M$ and consider the de Rham complex
$$
\Omega^*:=\big(\hdots \to 0 \to 0 \to \Omega^0(M)\to \Omega^1(M)\to \hdots \Omega^k(M) \hdots \big)\;,
$$
where the nonzero terms are concentrated in negative degrees. Then $H$ takes $\Omega^*$ to the spectrum
$$H(\Omega^*(M))=\left\{
\begin{array}{l}
DK\left(\hdots \to 0 \to 0 \to \Omega_{\rm cl}^0(M)\right)
\\
\\
DK\left(\hdots \to 0 \to 0 \to \Omega^0(M)\to \Omega_{\rm cl}^1(M)\hdots\right)
\\
\\
DK\left(\hdots \to 0 \to 0 \to \Omega^0(M)\to \Omega^1(M)\to \Omega^2_{\rm cl}(M) \hdots\right)
\\
\vdots
\\
DK\left(\hdots \to 0 \to 0 \to \Omega^0(M)\to \Omega^1(M)\to \hdots \to \Omega_{\rm cl}^k(M) \hdots\right) 
\\
\vdots
\end{array}\right.
$$
By the basic properties of the Dold-Kan functor, the stable homotopy groups of this spectrum are computed as
\bea
\pi_n^sH(\Omega^*(M)) &\simeq & \lim_{k\to \infty}\pi_{k+n}DK\big(\hdots \to 0 \to 0 \to \Omega^0(M)\to \Omega^1(M)\to \hdots \Omega_{\rm cl}^k(M)\big) 
\\
&\simeq &\lim_{k\to \infty}H_{k+n}\big(\hdots \to 0 \to 0 \to \Omega^0(M)\to \Omega^1(M)\to \hdots \Omega_{\rm cl}^k(M)\big)\;.
\eea
For $n>0$, these groups are $0$. For $n \leq 0$, they are the $n$th de Rham groups $H_{\rm dR}^n(M)$.
\end{example}

Now the functor $H$ in \eqref{EM fun} prolongs to a functor on prestacks
$$H:[\cartsp,\mathscr{C}{\rm h}]\to [\cartsp,\mathscr{S}{\rm p}]\;.$$
In fact, using the properties of the Dold-Kan correspondence, it is fairly straightforward to show that this functor preserves local weak equivalences \cite{Br}. 
We therefore get a functor of smooth stacks
\(
H:Sh_{\infty}(\cartsp;\mathscr{C}{\rm h})\to Sh_{\infty}(\cartsp;\mathscr{S}{\rm p})\;.
\label{Smooth EM}
\)
Recall that for an $\Omega$-spectrum $\E$, we always have a rational equivalence:
$$
{\rm r}:\E\wedge M\RR\to H\left(\pi_*(\E)\otimes \RR\right)\;,
$$
where $M\RR$ denotes an Eilenberg-Moore spectrum. Now, since we are working over the site of Cartesian spaces,
 the Poincar\'e lemma implies that the inclusion 
$j:\RR[0]\to \Omega^*$
induces an equivalence
$$
{\rm id} \otimes j:\pi_*(\E)\otimes \RR[0]\to \pi_*(\E)\otimes \Omega^*\;,
$$
where $\pi_*(\E)=\E(*)$ (which follows from suspension). 

\begin{definition}
\label{differential function}
Let $\E$ be a spectrum. For an unbounded chain complex $C_{\bullet}$, let $\tau_{\leq 0}C_{\bullet}$ 
denote the truncated complex 
$$
\tau_{\leq 0}C_{\bullet}=\left(\hdots 0 \to 0 \to 0 \to C_0\to C_{-1} \to \hdots \to C_{-n} \to \hdots \right)\;.
$$
A differential function spectrum ${\rm diff}(\E,{\rm ch})$ is a homotopy pullback 
$$\label{pullback diff}
\xymatrix{
{\rm diff}(\E,{\rm ch}) \ar[r]\ar[d] & H\left(\tau_{\leq 0} \Omega^*\otimes \pi_*(\E)\right)\ar[d]
\\
{\E} \ar[r]^-{\rm ch} & H\left(\Omega^*\otimes \pi_*(\E)\right)\;,
}
$$
where ${\rm ch}=j\circ {\rm r}$ and $j$ induces an equivalence
$
j:\pi_*(\E)\otimes \RR[0]\overset{\simeq}{\longrightarrow} \pi_*(\E)\otimes \Omega^*
$.
\end{definition}

\begin{remark} In our definition, we have chosen the complex $\Omega^*\otimes \pi_*(\E)$ as the de Rham complex modeling our rational theory. In general the differential function spectrum depends on this choice and on the equivalence $j$ \cite{Bun}. For the purposes of clarity and utility, we will always choose this model, although other models can be treated analogously. We do, however, keep the dependence on the map ${\rm ch}$ explicit to emphasize this fact. 
\end{remark}

\begin{example}[Deligne cohomology] 
Let $\E=H(\ZZ[n])\simeq \Sigma^nH\ZZ$ be the $n$-fold suspension of the Eilenberg-MacLane 
spectrum. In unbounded chain complexes, we have a natural isomorphism
$$
\underline{\ZZ}[n]\otimes \Omega^*\simeq \Omega^*[n]\;,
$$
where 
$\underline{\ZZ}[n]$ is the sheaf of locally constant integer-valued functions in degree $n$,
and the complex on the right hand side
 has been shifted up $n$ units. That is $\Omega^n$ is in degree $0$, 
while $\Omega^0$ is in degree $n$. 
Since $\Sigma^nH\ZZ$ is in the image of the Eilenberg-MacLane functor $H$ and $H$ 
preserves homotopy pullbacks, the homotopy pullback 
$$
\xymatrix{
{\rm diff}(\Sigma^nH\ZZ,{\rm ch}) \ar[r] \ar[d] & H(\tau_{\leq 0}\Omega^*[n])\ar[d]
\\
\Sigma^nH\underline{\ZZ} \ar[r]^{\rm ch} & H(\Omega^*[n])
}
$$
is presented by the homotopy pullback of unbounded chain complexes
$$
\xymatrix{
\underline{\ZZ}[n] \times^h_{\Omega^*[n]}\tau_{\geq 0}\Omega^*[n] \ar[d] \ar[r]& \tau_{\leq 0}\Omega^*[n]\ar[d]
\\
\underline{\ZZ}[n] \ar[r] & \Omega^*[n]\;.
}
$$
By stability, we can identify the homotopy pullback with the shifted mapping cone
$$
\underline{\ZZ}[n] \times^h_{\Omega^*[n]}\tau_{\leq 0}\Omega^*[n]\simeq {\rm cone}\left(\underline{\ZZ}[n]
\oplus \tau_{\leq 0}\Omega^* \to  \Omega^*[n]\right)[-1]\;.
$$
The right hand side is precisely the Deligne complex $\ZZ^{\infty}_{\cal D}(n+1)$. We therefore have an equivalence
$$
H(\ZZ^{\infty}_{\cal D}(n+1))\simeq {\rm diff}(\Sigma^nH\ZZ,{\rm ch})\;.
$$
The underlying theory that this spectrum represents is precisely Deligne cohomology. In fact, by the Dold-Kan
 correspondence, we have an isomorphism of graded abelian groups
$$
\pi_0\hom_{\mathscr{C}{\rm h}}(N(C(\{U_i\}), \ZZ^{\infty}_{\cal D}(n+1))\simeq 
\pi_0\map(\Sigma^\infty M,{\rm diff}(\Sigma^n H\ZZ,{\rm ch}))\;.
$$
Here $N$ denotes the normalized Moore complex (adjoint to the Dold-Kan functor $DK$) and $C(\{U_i\})$ 
denotes the {\v C}ech nerve of some good open cover of $X$. The right hand side is simply the definition of 
${\rm diff}(\Sigma^n H\ZZ,{\rm ch})^0(M)$, while the left hand side is the shifted total complex of the {\v C}ech 
Deligne double complex. It therefore computes the degree $n$ Deligne cohomology $H^n(M;\ZZ^{\infty}_{\cal D}(n+1))$. 
\end{example}

The above example illustrates what exactly differential function spectra have to 
do with differential cohomology theories. The following definition can be found in \cite{BNV}.

\begin{definition}\label{diff coh}
Let $\E$ be a spectrum and let 
$${\rm ch}:\E\to H(\tau_{\leq 0}\Omega^*\otimes \pi_*(\E))\;,$$
be the Chern character map as in Definition \ref{differential function}.
The {\rm differential $\E$-cohomology of a manifold} is the smooth cohomology 
theory with degree $n$ component
$$\widehat{\E}^n(M)\simeq {\rm diff}(\Sigma^n \E,{\rm ch})^0(M)\;.$$
\end{definition}

Since, for each $n$, ${\rm diff}(\Sigma^n\E,{\rm ch})$ is a smooth spectrum it fits into a diamond diagram of the form \eqref{spectra diamond}, as established in \cite{BNV}\cite{Sch}. In \cite{BNV}, it was shown that the form that this diamond takes is precisely the differential cohomology diamond in Remark \ref{diff diamond}. In particular, Proposition \ref{diff-cohomology-cal} 
allows us to calculate the ${\rm diff}(\Sigma^n\E,{\rm ch})$ cohomology in degrees away from $0$ as
$$
{\rm diff}(\Sigma^n\E,{\rm ch})^q(M)=\left\{\begin{array}{ccc}
\E^{n+q}(M) && q>0,
\\
\\
\E_{U(1)}^{n-1+q}(M)&& q<0.
\end{array}\right.
$$

\section{The smooth Atiyah-Hirzebruch spectral sequence (AHSS)}
\label{Sec AHSS}

In this section, we describe general machinery to construct an Atiyah-Hirzebruch spectral sequence (AHSS) 
 from a smooth spectrum ${ \E}$. We also describe how to compare this spectral sequence to the 
classical AHSS spectral sequence for the underlying  theory ${ \Pi}\E$, in nice cases.

\subsection{Construction of the spectral sequence via {\v C}ech resolutions}
\label{Sec Cech}

The trick to describing the spectral sequence is to choose the right filtration on a fixed manifold. 
In the local (projective) model structure on smooth stacks, a natural choice arises: namely the {\it {\v C}ech-type filtration on good open covers}. This is indeed the most natural choice, since the maps which are weakly inverted in the local model structure are precisely those arising from taking the {\v C}ech nerve of a good open cover of a manifold. That is, we have a weak equivalence
$$
w:{\rm hocolim}\left\{
\xymatrix{
\hdots \ar@<.75pc>[r] \ar@<.25pc>[r] \ar@<-.25pc>[r] \ar@<-.75pc>[r]  & \ar@<.5pc>[l] \ar@<-.5pc>[l] \ar[l]   
\coprod_{\alpha\beta\gamma}U_{\alpha\beta\gamma} \ar@<.5pc>[r] \ar@<-.5pc>[r] \ar[r] & \ar@<.25pc>[l] \ar@<-.25pc>[l] \coprod_{\alpha\beta}U_{\alpha\beta} \ar@<.25pc>[r] \ar@<-.25pc>[r] & \ar[l]\coprod_{\alpha}U_{\alpha}}
\right\}\to X\;.
$$
We now explicitly describe a filtration on $C(\{U_i\})$. 
Recall that any simplicial diagram $J:\Delta^{\rm op}\to Sh_{\infty}(\cartsp)$ can be filtrated by skeleta. 
More precisely, let $i:\Delta_{\leq k}\into \Delta$ denote the embedding of the full subcategory of linearly 
ordered sets $[r]$, such that $r\leq k$. Then $i$ induces a restriction between functor categories (the $k$-th truncation)
$$
\tau_{\leq k}:[\Delta^{\rm op},Sh_{\infty}(\cartsp)]\longrightarrow [\Delta_{\leq k}^{\rm op},Sh_{\infty}(\cartsp)]\;.
$$
By general abstract nonsense (the existence of left and right Kan extensions), there are left and right adjoints 
$({\rm sk}_k \dashv \tau_{\leq k} \dashv {\rm cosk}_k)$
$$
\xymatrix{
[\Delta^{\rm op},Sh_{\infty}(\cartsp)]\ar[rr]^{\tau_{\leq k}}~ &&~
[\Delta_{\leq k}^{\rm op},Sh_{\infty}(\cartsp)] \ar@<1em>[ll]_{{\rm cosk}_k}\ar@<-1em>_{{\rm sk}_k}[ll]
}.$$
Furthermore, by composing adjoints, we have an adjunction $({\rm sk}_k\dashv {\rm cosk}_k)$
$$
\xymatrix{
[\Delta^{\rm op},Sh_{\infty}(\cartsp)]\ar@<.5em>[rr]^{{\rm sk}_{k}} ~&&~ 
[\Delta^{\rm op},Sh_{\infty}(\cartsp)] \ar@<.5em>[ll]_{{\rm cosk}_k}\;.
}$$
The functor ${\rm sk}_k$ freely fills in degenerate simplices above level $k$, while ${\rm cosk}_k$
 probes a simplicial object with simplices only up to level $k$ (the singular $k$-skeleton).

\begin{proposition}
Let ${Y}_{\bullet}$ be a simplicial object in $Sh_{\infty}(\cartsp)$. Then we can filter ${Y}_{\bullet}$ by skeleta
$${\rm sk}_0{Y}_{\bullet} \to {\rm sk}_1{Y}_{\bullet} \to \hdots {\rm sk}_k{Y}_{\bullet} \to \hdots \to {Y}_{\bullet}\;.$$
The homotopy colimit over ${Y}_{\bullet}$ is presented by the ordinary colimit
$$
\underset{\Delta^{op}}{\rm hocolim}({Y}_{\bullet})\simeq 
\underset{k\to \infty}{\rm colim}~\underset{\Delta^{op}}{{\mathbb L}{\rm colim}}({\rm sk}_k{Y}_{\bullet})\;,
$$
where ${\mathbb L}{\rm colim}$ is the left derived functor of the colimit, hence computable upon suitable cofibrant replacement of the diagram \footnote{We take this particular model of the homotopy colimit in order to ensure that taking the colimit of the resulting diagram makes sense. The claim will also hold for other presentations of the homotopy colimit}. 
\end{proposition} 
\theproof
Since $Sh_{\infty}(\cartsp)$ is presented by a combinatorial simplicial model category, the homotopy colimit over a filtered diagram is presented by the ordinary colimit and the canonical map
$$\underset{k\to \infty}{{\mathbb L}{\rm colim}} ~\underset{\Delta^{\rm op}}{\mathbb{L}{\rm colim}}({\rm sk}_k{Y}_{\bullet})\to \underset{k\to \infty}{\rm colim} ~\underset{\Delta^{\rm op}}{\mathbb{L}{\rm colim}}({\rm sk}_k{Y}_{\bullet})$$
is an equivalence. 
Since homotopy colimits commute with homotopy colimits, we also have an equivalence
$$
\underset{k\to \infty}{{\mathbb L}{\rm colim}} ~\underset{\Delta^{\rm op}}{\mathbb{L}{\rm colim}}({\rm sk}_k{Y}_{\bullet})\simeq \underset{\Delta^{\rm op}}{\mathbb{L}{\rm colim}}~\underset{k\to \infty}{{\mathbb L}{\rm colim}}({\rm sk}_k{Y}_{\bullet})\;.
$$
Again, using the fact that the ordinary colimit over a filtered diagram presents the homotopy colimit, we have an equivalence
$$
\underset{\Delta^{\rm op}}{\mathbb{L}{\rm colim}}~\underset{k\to \infty}{{\mathbb L}{\rm colim}}({\rm sk}_k{Y}_{\bullet})\to \underset{\Delta^{\rm op}}{\mathbb{L}{\rm colim}}~\underset{k\to \infty}{{\rm colim}}({\rm sk}_k{Y}_{\bullet})\simeq \underset{\Delta^{\rm op}}{\mathbb{L}{\rm colim}}({Y}_{\bullet})\;.
$$
\vspace{-4mm}
\endofproof

\begin{remark}
 The above proposition says that the homotopy colimit over the simplicial object is filtered by homotopy colimits of its skeleta. In particular, if $M$ is a paracompact manifold, we can fix a good open cover on $M$ and form the simplicial object given by its {\v C}ech nerve 
$$ 
C(\{U_i\}):=
\xymatrix{
\hdots \ar@<.75pc>[r] \ar@<.25pc>[r] \ar@<-.25pc>[r] \ar@<-.75pc>[r]  & \ar@<.5pc>[l] \ar@<-.5pc>[l] \ar[l]   ~\coprod_{\alpha\beta\gamma}U_{\alpha\beta\gamma} \ar@<.5pc>[r] \ar@<-.5pc>[r] \ar[r] & \ar@<.25pc>[l] \ar@<-.25pc>[l] ~\coprod_{\alpha\beta}U_{\alpha\beta} \ar@<.25pc>[r] \ar@<-.25pc>[r] & \ar[l] ~\coprod_{\alpha}U_{\alpha}}
\;.$$
The homotopy colimit over this object is then filtered by its skeleta. 
\end{remark}

Let us see exactly what the skeleta look like in this case. To this end, we recall that in $Sh_{\infty}(\cartsp)$ the full homotopy colimit is presented by the local homotopy formula
$$
{\rm hocolim}_{\Delta^{op}}\ C(\{U_i\})=\int^{n\in \Delta}\coprod_{\alpha_0\hdots \alpha_n}U_{\alpha_0\hdots \alpha_n}\odot \Delta[n]\;.
$$ 
The filtration on this object is given by first truncating the {\v C}ech nerve and then freely filling in degenerate simplices. As a consequence, in degree $k$ we can forget about the simplices of dimension higher than $k$. The homotopy colimit over this skeleton is then given by a strict colimit over the diagram
{
\small 
\(
\hspace{-.6cm}
\label{k skeleton}
\xymatrix{
\coprod_{\alpha_0\hdots \alpha_k}U_{\alpha_0\hdots \alpha_k}\odot \Delta[k] \hdots \ar@<.75pc>[r] \ar@<.25pc>[r] \ar@<-.25pc>[r] \ar@<-.75pc>[r]  & \ar@<.5pc>[l] \ar@<-.5pc>[l] \ar[l]   \coprod_{\alpha\beta\gamma}U_{\alpha\beta\gamma}\odot \Delta[2] \ar@<.5pc>[r] \ar@<-.5pc>[r] \ar[r] & \ar@<.25pc>[l] \ar@<-.25pc>[l] \coprod_{\alpha\beta}U_{\alpha\beta}\odot \Delta[1] \ar@<.25pc>[r] \ar@<-.25pc>[r] & \ar[l]\coprod_{\alpha}U_{\alpha}}\odot \Delta[0]
\;,\)}
where the face and degeneracy maps are induced by the face and degeneracy maps of $\Delta[k]$. Taking $k\to \infty$, we do indeed reproduce the coend representing the full homotopy colimit $C(\{U_i\})$. 

We would like to eventually use this filtration to define a Mayer-Vietoris like spectral sequence for general cohomology theory 
${\E}$. To get to this step, however, we will need to identify the successive quotients of the filtration. To simplify notation in what follows, we will fix a manifold $M$ with {\v C}ech nerve $C(\{U_i\})$ and we set 
$$
X_k:=\underset{\Delta^{op}}{\rm hocolim}\big({\rm sk}_kC(\{U_i\})\big)\;.
$$
Then the quotient $X_k/X_{k-1}$ can be identified from the previous discussion by quotienting out the face maps at level $k$ described in diagram \eqref{k skeleton}. Since the tensor of a simplicial set and a stack is given by the product of the stack with the discrete inclusion of the simplicial set, we can identify the quotient from the pushout of coends
$$\xymatrix{
\int^{n< k}\underset{\alpha_0\hdots \alpha_{n}}{\coprod} U_{\alpha_0\hdots \alpha_{n}}\times {\rm disc}(\Delta[n])\ar[d]^{\partial}\ar[r] & \ast
\\
\int^{m\leq k}\underset{\alpha_0\hdots \alpha_m}{\coprod}U_{\alpha_0\hdots \alpha_m}\times {\rm disc}(\Delta[m])\;,&
}
$$
where $\partial$ denotes the boundary inclusion. At the level of points (or elements), a simplex in $\int^{n< k}\coprod_{\alpha_0\hdots \alpha_{k}}U_{\alpha_0\hdots \alpha_{n}}\times {\rm disc}(\Delta[n])$ is given by a pair 
$$
(\rho,\sigma)\in \coprod_{\alpha_0\hdots \alpha_{k-1}}U_{\alpha_0\hdots \alpha_{k-1}}\times {\rm disc}(\Delta[k-1])\;,
$$ 
which is glued to lower simplices via the face and degeneracy relations. 

Let us identify where the boundary inclusion takes a generic simplex. Then the quotient $X_k/X_{k-1}$ will be obtained by gluing these simplices together to a single point. Note that the face and degeneracy relations imply that simplices of the form $(\rho, s_{j+1}\sigma)$ are sent by $d_j$ to $(d_j\rho, \sigma)$. Since simplices in the image of the face maps are precisely those which are collapsed to a point, we see that 
$$
(d_j\rho, \sigma)\sim \ast\ \ \ \ {\rm for~every~} \sigma.
$$ 
We therefore see that each term of the coproduct $\coprod_{\alpha_0\hdots \alpha_k}U_{\alpha_0\hdots \alpha_k}$ is joined to another by the inclusion into a lower intersection. These lower intersections are then collapsed to a point yielding the wedge product
$$
\bigvee_{\alpha_0\hdots \alpha_k}U_{\alpha_0\hdots \alpha_k}\subset X_k/X_{k-1}\;.
$$
Similarly, the simplex $(s_{j+1}\rho,\sigma)$ is sent to $(\rho,d_j\sigma)$ under $d_j$. We therefore identify the discrete simplicial sphere in the quotient
$$
{\rm disc}(\Delta[k]/\partial \Delta[k])\subset X_k/X_{k-1}\;.
$$
Finally, the relations imposed by the coend imply that a simplex of the form $(s_{j}\rho,\sigma)$ is glued to $(\rho, d_{j}\sigma)$. The former are precisely those simplices in the simplicial sphere while the later are glued to the point. Similarly, 
$(\rho, s_j\sigma)$ is glued to the point. Thus we have the following.

\begin{lemma}
\label{Lemma vee}
We can identify the quotient with the smash product
$$
X_k/X_{k-1}\simeq {\rm disc}(\Delta[k]/\partial \Delta[k])\wedge \bigvee_{\alpha_0\hdots \alpha_k}U_{\alpha_0\hdots \alpha_k}\simeq \Sigma^k\Big(\bigvee_{\alpha_0\hdots \alpha_k}U_{\alpha_0\hdots \alpha_k}\Big)\;.
$$
\end{lemma}

\begin{remark} [The filtration as a natural choice] 
Another way to think of our filtration above is the following. 
  Let us form a {\v C}ech nerve of a manifold, then contract
all the patches and intersections in that {\v C}ech nerve as points, such
as to just obtain a simplicial set. Then the {\it Borsuk's nerve theorem} 
(see \cite{Bj} for a survey, \cite{Ha} Corollary 4G.3, or \cite{Pr} Theorem 3.21)
 says that this simplicial set is
equivalent -- weak homotopy equivalent -- to the singular simplicial
complex of the manifold, hence to its homotopy type.
Moreover, that singular simplicial complex (or rather
its geometric realization), in turn, gives a CW-complex realization of the
original manifold.
So with this in mind, one may view our  filtration above 
 as the natural smooth refinement of the filteration by
CW-stages of the manifold. That is, in taking the {\v C}ech nerve
{\rm without} contracting all its patches to points, we
retain exactly the smooth information that, via Borsuk's theorem,
corresponds to each cell in the canonical CW-complex incarnation of
the manifold. So in this sense, our  refinement can be viewed as the 
canonical smooth refinement of the traditional filtering by CW-stages.

\label{Rem Borsuk}
\end{remark}

We are now ready to describe the spectral sequence.

\begin{theorem}[AHSS for general smooth spectra]\label{main theorem}
Let $M$ be a compact smooth manifold and let ${\E}$ be a smooth spectrum. There is a spectral sequence with 
$$
E_{2}^{p,q}=H^p(M,{ \E}^{q})\ \Longrightarrow \ { \E}^{p+q}(M)\;.
$$
Here $H^p$ denotes the $p$-th {\v C}ech cohomology with coefficients in the presheaf ${ \E}^{q}$. Moreover, the differential on the $E_1$-page is given by the differential in {\v C}ech cohomology.
\label{MV thm} 
\end{theorem} 
\theproof
The proof is almost immediate from the definitions. Recall that we have identified the quotients
in Lemma \ref{Lemma vee}.
By the axioms for a smooth cohomology theory, we have that the ${\E}$-cohomology of the quotient is given by
\bea
{ \E}^*(X_k/X_{k-1})&\simeq& { \E}^*\Big(\Sigma^k\Big(\bigvee_{\alpha_0\hdots \alpha_k}U_{\alpha_0\hdots \alpha_k}\Big)\Big)
\\
&\simeq & { \E}^{*-k}\Big(\bigvee_{\alpha_0\hdots \alpha_k}U_{\alpha_0\hdots \alpha_k}\Big)
\\
&\simeq &  \bigoplus_{{\alpha_0\hdots \alpha_k}}{\E}^{*-k}\Big(U_{\alpha_0\hdots \alpha_k}\Big)\;.
\eea
Applying ${ \E}^{p+q}$ to the cofiber squence
$
X_p\into X_{p+1}\onto X_{p+1}/X_p
$
gives the long exact sequence in ${ \E}$-cohomology
\(
\hdots { \E}^{p+q}(X_{p+1}/X_p) \to { \E}^{p+q}(X_{p+1})\to { \E}^{p+q}(X_p)\to { \E}^{p+q+1}(X_{p+1}/X_p) \hdots\;. 
\label{LES XX}
\)
Forming the corresponding exact triangle, we get a spectral sequence with $E^{p,q}_1$ term
$$
E^{p,q}_1=\bigoplus_{\alpha_0,\hdots, \alpha_p}{ \E}^{q}(U_{\alpha_0\hdots \alpha_p})\;.
$$
Now we want to show that the differential on this page is given by the {\v C}ech differential 
$$
\delta:E^{p,q}_1=\bigoplus_{\alpha_0\hdots \alpha_p}{ \E}^{q}(U_{\alpha_0\hdots \alpha_p})~\longrightarrow~\bigoplus_{\alpha_0\hdots \alpha_{p+1}}{\E}^{q}(U_{\alpha_0\hdots \alpha_{p+1}})=E^{p+1,q}\;.
$$
To this end, note that differential on the $E_1$-page, by definition, comes from the exact sequence
 $$\hdots \to \E^{p+q}(X_{p+1}/X_p)\overset{j}{\to} \E^{p+q}(X_{p+1})\overset{i}{\to} \E^{p+q}(X_p)\overset{\partial}\to \E^{p+q+1}(X_{p+1}/X_p)\to \hdots\;.$$
We need to show that $\partial j=d_1=\delta$ is the {\v C}ech differential. By naturality of the connecting homomorphism $\partial$, we have a commutative diagram
$$
\hspace{-.6cm}\label{cech differential}
\xymatrix{
\check{C}^{p-1}(M;\E^q)\ar[d]_{\simeq} \ar[rr]^{d_1} && \check{C}^{p}(M;\E^{q})\ar[d]^{\simeq}
\\
\bigoplus_{\alpha_0\hdots \alpha_{p-1}}\E^q(U_{\alpha_0\hdots \alpha_{p-1}}) \ar[d]_{\simeq}\ar[rr] && \bigoplus_{\alpha_0\hdots \alpha_{p}}\E^q(U_{\alpha_0\hdots \alpha_{p}})\ar[d]^{\simeq}
\\
 \E^{p+q-1}(X_{p-1}/X_{p-2}) \ar[r]^-{j}\ar[d] & \E^{p+q-1}(X_{p-1})\ar[r]^-{\partial}\ar[d] & \E^{p+q}(X_{p}/X_{p-1})\ar[d]
 \\
 \E^{p+q-1}(\partial\Delta[p]\times U_{\alpha_0\hdots \alpha_{p-1}}) \ar[r]^{\rm id} & \E^{p+q-1}(\partial\Delta[p]\times U_{\alpha_0\hdots \alpha_{p-1}})\ar[r]^-{\partial} & \E^{p+q}\Big(\Delta[p]/\partial\Delta[p]\wedge U_{\alpha_0\hdots \alpha_{p}}\Big)\;,
 }
 $$
  where the vertical bottom maps are induced from the inclusion of a factor
  \(\label{face maps}
  \xymatrix{
  \Delta[p]\times U_{\alpha_0\hdots \alpha_{p}}\ar@{^{(}->}[r] &  X_{p}
  \\
  \partial\Delta[p]\times U_{\alpha_0\hdots \alpha_{p-1}}\ar@{^{(}->}[r]\ar[u] & X_{p-1} \ar[u]
  \\
  \emptyset \ar[r]\ar[u] & X_{p-2}\ar[u]
  }
  \)
  into the $p$-level of the filtration. Comparing the top and bottom composite morphisms in the big diagram, we see that on $(p-1)$-fold intersections $U_{\alpha_0\hdots \alpha_{p-1}}$, the map $d_1$ is forced to map a section to the alternating sum of restrictions, as this is precisely the map induced by the boundary inclusion in \eqref{face maps}.
  
All that remains is the convergence. To establish that, we simply note that compactness implies that, for large values of $p$, we have an equivalence $X_p\simeq X$. Moreover, there are only finitely many diagonal entries at each page of the sequence. With this assumption, the convergence to the corresponding graded complex
$$
E_{\infty}^{p,q}=\frac{{\rm ker}\left(\E^{p+q}(X) \to \E^{p+q}(X_{p})\right)}{{\rm ker}\left(\E^{p+q}(X) \to \E^{p+q}(X_{p+1})\right)}=\frac{F_p \E^{p+q}(X)}{F_{p+1}\E^{p+q}(X)}
$$
follows exactly as in the classical case in \cite{AH}.  
\endofproof

\paragraph{Fiber bundles.}We can also construct a spectral sequence for a fiber bundle
$$
F\to N \overset{p}{\to} M\;,
$$
where each map is a smooth map of manifolds and $M$ is compact. To that end, we note that for a fixed good open cover $\{U_i\}$ of $M$, the pullbacks $\{p^{-1}(U_i)\}$ define a good open cover of $N$. By local triviality, we have that each $p^{-1}(U_i)\simeq F\times U_i$. Then, using the filtration
$$
X_k=\underset{\Delta^{\rm op}}{\rm hocolim}\big({\rm sk}_kC(\{p^{-1}(U_i)\})\big)
$$
on the total space $N$, we identify the successive quotients
$$
X_k/X_{k-1}\simeq \Sigma^k \bigvee_{\alpha_0\hdots \alpha_k} U_{\alpha_0\hdots \alpha_k}\wedge F\;.
$$
A similar argument as in the proof of Theorem \ref{MV thm} gives

\begin{theorem}[Smooth AHSS for fiber bundles]
Let $M, N$ and $F$ be manifolds, with $M$ compact. Let $F\to N \overset{p}{\to} M$ be a fiber bundle. Let ${ \E}$ be a sheaf of spectra. Then there is a spectral sequence
$$
E_{2}^{p,q}=H^p(M,{ \E}^{q}(-\wedge F ))\ \Longrightarrow \ { \E}^{p+q}(N)\;.
$$
Here $H^p$ denotes the $p$-th {\v C}ech cohomology with coefficients in the presheaf ${ \E}^{-q}(-\wedge F)\;.$
\end{theorem}

\begin{remark}[Unreduced theories] Note that the smooth spectral sequence works for reduced theories. One can treat unreduced theories similarly by setting 
$$
\E^q(M,\ast):=\tilde{\E}^q(M_{+})\;,
$$
where the tilde denotes the reduced theory and $M_+$ is the pointed stack with basepoint $\ast$. In this case, we have the slight modification on the second spectral sequence, which takes the form
$$
E_{2}^{p,q}=H^p(M,{ \E}^{q}(-\times F ))\ \Longrightarrow \ { \E}^{p+q}(N)\;.
$$
\end{remark}

\subsection{Morphisms of smooth spectral sequences and refinement of the AHSS}
\label{Sec Morph}

Our next task will be to show that these spectral sequences do indeed refine the classical 
Atiyah-Hirzebruch spectral sequence (AHSS) \cite{AH}. 
Since any smooth theory $\E$ comes as a refinement of the underlying topological theory ${ \Pi} \E$, 
we will immediately get a morphism of spectral sequences induced by the morphism of spectra
$$
I:\E\to { \Pi}\E\;.
$$
Unfortunately, this morphism does not allow us to compare the differentials of the spectral sequences in the 
way  that we would ideally hope for. However, as we will progressively see, the situation can be remedied by 
constructing a slightly different morphism of spectral sequences. This morphism is related to the 
\emph{boundary map} of spectral sequences which occurs when a morphism of spectra induces the $0$ 
map on corresponding spectral sequences (see \cite{Mi} for a discussion in the case of the Adams spectral sequence). 
We first discuss the morphism induced by $I$ and then construct this ``boundary type" map and prove that
 it indeed defines a morphism of spectral sequences.

\begin{definition}
Let $E^{p,q}_n$ and $F^{p,q}_n$ be spectral sequences, that is, a sequence of bigraded complexes 
$E^{p,q}_n$ and $F^{p,q}_n$, $n\in \NN$. A {\rm morphism of spectral sequences} is a morphism of 
bigraded complexes
$$
f_n:E^{p,q}_n\to F^{p,q}_n\;,
$$
defined for all $n>N$, where $N$ is some fixed positive integer. Furthermore, we require the map 
$f_{n+1}$ to be the map on homology induced by $f_n$. We call the smallest integer $N$ such that 
$f_n$ are defined for $n>N$ the \emph{rank} of the morphism. 
\end{definition}

We now apply this to the smooth AHSS. The next result should follow 
from general principles, but we emphasise it explicitly for clarity and for subsequent use. 

\begin{proposition} Let $\E$ and ${\cal F}$ be smooth spectra. Then a map $f:\E\to {\cal F}$ 
induces a morphism of corresponding smooth AHSS's
$$E^{p,q}_n\to F^{p,q}_n\;.$$
\end{proposition}
\theproof
Fix a manifold $X$ and a good open cover $\{U_i\}$. Let $X_p$ denote the $p$-th filtration of the {\v C}ech nerve as before. It is clear by naturality that a map of spectra $f:\E\to {\cal F}$ induces a morphism of long exact sequences (see 
\eqref{LES XX})
$$
\xymatrix{
\hdots { \E}^{p+q}(X_{p+1}/X_p) \ar[r]\ar[d] & { \E}^{p+q}(X_{p+1})\ar[r]\ar[d] &  { \E}^{p+q}(X_p) \ar[r]\ar[d] & { \E}^{p+q+1}(X_{p+1}/X_p)\ar[d] \hdots 
\\
\hdots {\cal F}^{p+q}(X_{p+1}/X_p) \ar[r] & {\cal F}^{p+q}(X_{p+1})\ar[r] &  {\cal F}^{p+q}(X_p) \ar[r] & {\cal F}^{p+q+1}(X_{p+1}/X_p) \hdots 
\;.}$$
It follows immediately from the construction of the corresponding exact triangles that this morphism commutes with the differentials.
\endofproof

This now allows us to compare the topological and the smooth theories.

\begin{corollary}\label{refinement}
Let $\E$ be a smooth spectrum and ${ \Pi}\E$ be the underlying topological theory. 
Let $E_n$ and $F_n$ denote the spectral sequences corresponding to 
$\E$ and  ${ \Pi}\E$, respectively. 
The natural map $I:\E \to { \Pi}\E$ induces a morphism of classical AHSS's
\footnote{Here we have an unfortunate conflict of notation. We are using the 
same symbols for the pages in the spectral sequences for both the classical
and the refined theories. We will aim to make the context explicit whenever 
a possible ambiguity arises.}
$$
I: E^{p,q}_{n}\to F^{p,q}_n\;.
$$
\end{corollary}

\begin{remark}
It is interesting to note that the smooth spectrum ${ \Pi} \E$ is, by definition, locally constant. 
From the discussion around \eqref{cohesion}, this means that we have isomorphism 
$$
{ \Pi} E^q(U)
\simeq \pi_{-q}\map(U,{ \Pi}\E)
\simeq \pi_{-q}\map(*, { \Pi}\E)
\simeq \pi_{-q}{\Pi}\E
\simeq {\Pi}\E^q(\ast)
$$
for every element of a good open cover (or higher intersection) $U$. 
This connects, via Borsuk's theorem mentioned in Remark \ref{Rem Borsuk} above, 
the ``smooth
AHSS for locally constant coefficients" with the classical AHSS: the
locally constant coefficients see each (contractible) patch as a point,
and hence by Borsuk's theorem they see our ``{\v C}ech filteration" to be the
classical CW-cell filteration.
\end{remark}

From the construction of our smooth AHSS, 
we immediately get that the spectral sequence associated to the smooth spectrum 
is a refinement of the classical topological AHSS.

\begin{corollary} \label{refinement 1}
The spectral sequence $F^{p,q}_n$ is precisely the AHSS for the cohomology theory $\Pi \E$.
\end{corollary}

We now would like to apply the above machinery to differential cohomology theories. In particular, 
we note that for a differential function spectrum ${\rm diff}(\E,{\rm ch})$, the natural map 
$$
I:{\rm diff}(\E,{\rm ch})\to \underline{\E}\;,
$$
which strips the differential theory of the differential data and maps to the bare underlying theory,
is precisely the map induced by the unit $I:{\rm id}\to { \Pi}$. In the above discussion, we observed 
that this map always induces a morphism of spectral sequences. Moreover, the target spectral sequence 
is exactly the AHSS for the underlying topological theory. One might hope to be able to use this map to 
compare the differentials in the refined theory with those differentials in the classical AHSS.

Unfortunately, this does not work in practice, as we will see when we discuss applications in Sec. \ref{Sec App}. 
The core issue is that the spectral sequence for the refined theory usually ends up shifted with respect to 
the classical AHSS. As a consequence, the nonzero terms in each sequence are interlaced with respect to one 
another and the map $I$ ends up killing all the nonzero terms. This, in turn, stems from the appearance of 
the Bockstein map  (which raises degree by 1) in the differential cohomology diagram. 

However, there is often a different map between the \emph{lower quadrants} of the the two spectral sequence corresponding
 to ${\rm diff}(\E,{\rm ch})$ and $\E$, which lowers the degree as to match the corresponding nonzero entries. This map is 
 related to the so-called \emph{boundary map} between spectral sequences studied in \cite{Mi}. The next proposition concerns this map and will be essential for comparing the differentials in the refined theory to those of the classical theory.

\begin{proposition}
\label{Bockstein}
{ \bf (i)} Let $\E$ be a spectrum such that $\pi_*(\E)$ is concentrated in degrees which are a multiple of some integer 
$n\geq 2$ (e.g. K-theory, Morava K-theory). Suppose, moreover, that $\pi_*(\E)$ is projective in those degrees. 
Then the sequence of spectra
$$
\E\to \E\wedge M\RR \to \E \wedge MU(1)\overset{\beta_{\E}}{\longrightarrow} \Sigma \E\;,
$$
induces a short exact sequence on coefficients
\(\label{bockstein1}
0\to \pi_*(\E) \to \pi_*(\E)\otimes \RR \to \pi_{*}(\E) \otimes U(1)\to 0\;.
\)
{\bf (ii)} Let $\beta$ denote the connecting homomorphism (i.e. the Bockstein) for the coefficient sequence \eqref{bockstein1}. Let $E^{p,q}_n$ denote
 the spectral sequence corresponding to $\Sigma^{-1}\E \wedge MU(1)$ and let $F^{p,q}_n$ denote the spectral
  sequence corresponding to $\E$. Then 
$$
\beta: E^{p,q}_n\to F^{p,q}_n
$$
induces a morphism of spectral sequences of rank $2$.
\end{proposition}
\theproof
 Consider the long Bockstein sequence 
$$
\hdots \E \overset{r}{\longrightarrow} \E \wedge M\RR \overset{e}{\longrightarrow} 
\E \wedge MU(1) \overset{\beta_{\E}}{\longrightarrow} \Sigma \E\hdots\;,
$$
induced by the cofiber sequence
$$
\mathbb{S}\to M\RR \to MU(1)\;.
$$
Fix a manifold $M$ and let $X_p$ denote the $p$-level of the {\v C}ech filtration. Now each spectrum in the above sequence has a long exact sequence induced be the cofiber sequences
$$
X_{p-1}\to X_p\to X_p/X_{p-1}\;.
$$
from which one builds the exact couple for for the corresponding spectral sequence. Using the properties of $\pi_*(\E)$ along with this sequence, we can fit the long exact sequences into a diagram
$$
\hspace{-6mm}
\xymatrix{
\check{C}^p(X;\pi_{-q-1}(\E)) \ar[r]^-{q^*}\ar[d]^{r} & \E^{p+q-1}(X_{p})\ar[r]^-{i^*}\ar[d]^{r} &  \E^{p+q-1}(X_{p-1}) \ar[r]^-{\partial}\ar[d]^{r} & 0\ar[d] 
\\
\check{C}^p(X;\pi_{-q-1}(\E_{\RR})) \ar[r]^-{q^*}\ar[d]^{e} & \E_{\RR}^{p+q-1}(X_{p})\ar[r]^-{i^*}\ar[d]^{e} &  \E_{\RR}^{p+q-1}(X_{p-1}) \ar[r]^-{\partial}\ar[d]^{e} & 0\ar[d] 
\\
 \check{C}^p(X;\pi_{-q-1}(\E_{U(1)})) \ar[r]^-{q^*}\ar[d]^{\rm \beta_{\E}} & \E_{U(1)}^{p+q}(X_{p})\ar[r]^-{i^*}\ar[d]^{\beta_{\E}} &  \E_{U(1)}^{p+q}(X_{p-1}) \ar[r]^-{\partial}\ar[d]^{\beta_{\E}} & 0 \ar[d] 
\\
0 \ar[r] & {\E}^{p+q}(X_{p+1})\ar[r] &  {\E}^{p+q}(X_p) \ar[r] & \check{C}^p(X;\pi_{-q+1}(\E))\;, 
}
$$ 
where both the rows and columns are part of exact sequences and $\check{C}^p(X;A)$ denotes the group of {\v C}ech $p$-cochains with coefficients in $A$.  Since everything commutes, this induces a correponding short exact sequence of $E_1$ pages. At each $(p,q)$-entry this sequence is given by
$$
0\to C^{p}(X;\pi_{-q}(\E))\to C^{p}(X;\pi_{-q}(\E)\otimes \RR) \to C^{p}(X;\pi_{-q}(\E)\otimes U(1)) \to 0\;.
$$ 
Since the differentials on the $E_1$ page are precisely the {\v C}ech differentials, the construction of the Bockstein map 
in {\v C}each cohomology will produce a map of $E_2$-pages
$$
\beta: H^{p}(X;\pi_{-q}(\E)\otimes U(1))\to H^{p+1}(X;\pi_{-q}(\E))\;.
$$
We need to show that this map commutes with the differential. Choose a representative $x$ of a class in $H^p(X;\pi_{-q}(\E)\otimes U(1))$. By definition, $y=\beta(x)$ is a class such that $r(y)=\delta(\overline{x})$, where $\overline{x}$ is such that $e(\overline{x})=x$.  Then 
$$
r(d_2y)=d_2 r(y)=d_2 \delta(\overline{x})\;,
$$
We want to show that there is a lift $z$ of $d_2 x$ such that $\delta(z)=d_2\delta(\overline{x})$. Indeed, if this is the case, then $d_2y$ represents $\beta(d_2x)$ and we are done. 

To construct $z$, recall that $d_2x$ is defined by first pulling back by the quotient $q$, which lies in the image of the map 
induced by the inclusion $i:X_p\into X_{p+1}$, and then applying the boundary to an element of the preimage. Let $w$ be 
such that
$$
i^*(w)=q^*(x)\;.
$$
Chasing the diagram
$$
\xymatrix{
\check{C}^p(X;\pi_{-q-1}(\E_{\RR})) \ar[r]^-{q^*}\ar[d]^-{e} & \E_{\RR}^{p+q-1}(X_{p})\ar[d]^-{e} \ar[r]^-{e} & \E_{\RR}^{p+q-1}(X_{p-1})\ar[d]^-{e}
\\
\check{C}^p(X;\pi_{-q-1}(\E_{U(1)})) \ar[r]^-{q^*}\ar[d]^{\beta_{\E}} &  \E_{U(1)}^{p+q}(X_{p})\ar[d]^-{\beta_{\E}}\ar[r]^-{i^*} &  \E_{U(1)}^{p+q}(X_{p-1})\ar[d]^{\beta_{\E}}
 \\
 0 \ar[r] & \E^{p+q}(X_{p})\ar[r]^{i^*} & \E^{p+q}(X_{p-1})\;,
 }
$$ 
we see that $0=\beta_{\E}q^*(x)=\beta_{\E}i^*(w)=i^*(\beta_{\E} w)$. By exactness of the rows, this implies that $\beta_{\E}w=0$. Therefore, there is a class $\overline{w}\in \E_{\RR}^{p+q+1}(X_{p+1})$ such that $e(\overline{w})=w$. 

Now, by definition of the differential, we have
$$
e(\partial \overline{w})=\partial (e(\overline{w}))=\partial w=d_2x
$$
and $z:=\partial \overline{w}$ is a lift of $d_2x$. Using the fact that $\delta=d_1=\partial q^*$, we have
$$
\delta(z)=\delta(\partial \overline{w})=\partial (q^* \partial \overline{w})\;.
$$
By exactness, we have 
$$
i^*(q^*\partial \overline{w})=0=q^*\partial q^*(\overline{x})=q^*(\delta(\overline{x})),
$$
and it follows from the definition that $\delta(z)=d_2(\delta(\overline{x}))$.

To show that $H^*(\beta)$ commutes with the higher differentials, we proceed by induction. The above proves the base case. Suppose $\beta$ induces a map $H_n(\beta)$ on $E_n$ which commutes with $d_n$. Then $H^n(\beta)$ induces a well defined map $H^{n+1}(\beta)$ on the $E_{n+1}$ page. Let $x\in \bigcap _{i=1}^n\ker(d_{n+1})$ be a representative of a class on the $E_n$ page. Then be definition, $H^{n+1}(\beta)(x)=\beta(x)$ and the exact same argument as before (replacing $d_{2}$ with $d_{n+1}$), gives the result.
\endofproof

Having done the heavy lifting in the above proposition, we will now apply this to straightforwardly 
relate the differentials of the refined theory to those of the underlying topological theory. 
This will use an explicit alternative to the map $I$, along the lines of the discussion just before 
the statement of the above proposition.

\begin{theorem}[Refinement of differentials] 
\label{Thm beta}
Let $\E$ be a spectrum satisfying the properties of Proposition \ref{Bockstein} and let ${\rm diff}(\E,{\rm ch})$ 
be a differential function spectrum refining $\E$. Let $E_n$ and $F_n$ denote the smooth AHSS's corresponding 
to ${\rm diff}(\E,{\rm ch})$ and $\E$, respectively. Then the Bockstein $\beta$ defines a rank 2 morphism of 
fourth quadrant spectral sequences
$$
\beta: E^{p,q}_n\to F^{p,q}_n,\ \ q<0\;.
$$
\end{theorem}
\theproof
Recall that for $q<0$, Proposition 4 implies that ${\rm diff}(\E,{\rm ch})^q(M)\simeq \E^{q-1}_{U(1)}(M)$. The claim then follows from the previous proposition.
\endofproof

\subsection{Product structure and the differentials}
\label{Sec Product}

Let $\E$ be an $E_{\infty}$ ring spectrum. Then the associative graded-commutative product on $\E^*$ induces a product (associative and graded-commutative) on the refinement ${\rm diff}(\Sigma^n\E,{\rm ch})^*$,  
that is, a map
\(
\label{cup product}
\cup:~{\rm diff}(\Sigma^n\E,{\rm ch})^k\otimes {\rm diff}(\Sigma^m\E,{\rm ch})^j
\longrightarrow {\rm diff}(\Sigma^{n+m}\E,{\rm ch})^{k+j}
\)
(see \cite{Bun} \cite{U1}).
The goal of this section will be to establish the following very useful property, 
in analogy with the classical case.
\begin{proposition}[Compatibility with products]
\label{cup proposition}
The product 
$$
\cup:{\rm diff}(\Sigma^n\E,{\rm ch})^k\otimes {\rm diff}(\Sigma^m\E,{\rm ch})^j\to {\rm diff}(\Sigma^{n+m}\E,{\rm ch})^{k+j}
$$
induces a morphism of spectral sequences
$$
\cup:E_*(n) \times E_*(m)\to E_*(n+m)\;.
$$
Moreover, the differentials satisfy the Liebniz rule
$$
d(xy)=d(x)y+(-1)^{p+q}xd(y)\;.
$$
\end{proposition}

Let us first work out what the cup product pairing is on the $E_1$-page. Recall from the construction of the spectral sequence
 that the $E^{p,q}_1$ is given by
$$
E^{p,q}_1=\bigoplus_{\alpha_0\hdots \alpha_p}{\rm diff}(\Sigma^n\E,{\rm ch})^q(U_{\alpha_0\hdots \alpha_p})
\simeq \check{C}^p(M;{\rm diff}(\Sigma^n\E,{\rm ch})^q)\;.
$$
Using the product \eqref{cup product}, we get a cross product map
\begin{alignat}{2}
\hspace{-8mm}\times:\bigoplus_{\alpha_0\hdots \alpha_p}{\rm diff}(\Sigma^n\E,{\rm ch})^q(U_{\alpha_0\hdots \alpha_p})
&\times \bigoplus_{\alpha_0\hdots \alpha_r}{\rm diff}(\Sigma^m\E,{\rm ch})^t(U_{\alpha_0\hdots \alpha_r})\to
\nonumber\\
&\to
 \bigoplus_{\alpha_0\hdots \alpha_p}\bigoplus_{\alpha_0\hdots \alpha_r}{\rm diff}(\Sigma^{n+m}\E,{\rm ch})^{q+t}(U_{\alpha_0\hdots \alpha_p}\times U_{\alpha_0\hdots \alpha_r})\;.
\end{alignat}
We also have an isomorphism
\begin{alignat}{2}
\hspace{-3mm}
\bigoplus_{\alpha_0\hdots \alpha_s}{\rm diff}(\Sigma^{n+m}\E,{\rm ch})^{q+t}((U\times U)_{\alpha_0\hdots \alpha_{s}}) &\simeq  {\rm diff}(\Sigma^{n+m}\E,{\rm ch})^{q+t}\Big(\bigvee_{\alpha_0\hdots \alpha_s}(U\times U)_{\alpha_0\hdots \alpha_{s}}\Big) 
\nonumber\\
\nonumber\\
&\simeq {\rm diff}(\Sigma^{n+m}\E,{\rm ch})^{q+t}\Big(\bigvee_{\alpha_0\hdots \alpha_p}\bigvee_{\alpha_0\hdots \alpha_r}\bigvee_{p+r=s}U_{\alpha_0\hdots \alpha_p}\times U_{\alpha_0\hdots \alpha_r}\Big)
\nonumber\\
\nonumber\\
&\simeq \bigoplus_{\alpha_0\hdots \alpha_p}\bigoplus_{\alpha_0\hdots \alpha_r}\bigoplus_{p+r=s}{\rm diff}(\Sigma^{n+m}\E,{\rm ch})^{q+t}(U_{\alpha_0\hdots \alpha_p}\times U_{\alpha_0\hdots \alpha_r})\;,
\nonumber
\end{alignat}
given by decomposing the product of the cover $\{U_{\alpha}\}$ with itself. Finally, we can pullback by the diagonal map
\bea
\Delta^*:\bigoplus_{\alpha_0\hdots \alpha_s}{\rm diff}(\Sigma^{n+m}\E,{\rm ch})^{q+t}((U\times U)_{\alpha_0\hdots \alpha_{s}})\to \bigoplus_{\alpha_0\hdots \alpha_s}{\rm diff}(\Sigma^{n+m}\E,{\rm ch})^{q+t}(U_{\alpha_0\hdots \alpha_{s}})\simeq
\\
\simeq
 \check{C}^{p+r}(M;{\rm diff}(\Sigma^{n+m}\E,{\rm ch})^{q+t})\;.
\eea
The cup product on the $E_1$-page is defined by the composite map $\Delta^*\times$.
\begin{lemma}
\label{lem d1}
The differential $d_1$ on the $E_1$-page satisfies the Leibniz rule.
\end{lemma}
\theproof
The construction of the cup product on the $E_1$-page is precisely the cup product structure 
for {\v C}ech-cohomology. The {\v C}ech differential satisfies the Leibniz rule and this is precisely 
$d_1$, by construction.
\endofproof

We are now ready to prove Proposition \ref{cup proposition}.

\theproof
The proof follows by induction on the pages of the spectral sequence. The base case is satisfied by 
Lemma \ref{lem d1}. Now suppose we have a cup product map 
$$
\cup: E(n)_k\times E(n)_k\to E(n+m)_k\;,
$$ 
such that $d_k$ satisfies Leibniz. By definition, we have
$$
E(n)^{p,q}_{k+1}=\frac{{\ker\big(d_k:E(n)^{p,q}_k\to E(n)^{p+k,q+k-1}_k\big)}}
{{\rm im}\big(d_{k}:E(n)^{p-k,q-k+1}\to E(n)^{p,q}\big)}\;,
$$
and we define the cup product
$$
\cup:E(n)^{p,q}_{k+1}\times E(m)^{r,s}_{k+1}\to E(n+m)^{p+r,q+s}_{k+1}
$$
by restricting to elements in the kernel of $d_k$. The product is well defined since $d_k$ satisfies the Leibniz rule. At this stage the problem looks formally like the classical problem. 
Hence, analogously to the classical discussion in \cite{Ha}, it is tedious but straightforward to 
show that $d_{k+1}$ also satisfies the Leibniz rule.
\endofproof

\section{Applications to differential cohomology theories}
\label{Sec App}

In this section we  would like  to apply the spectral sequence constructed in the previous section 
to various differential  cohomology theories. The construction is general enough to apply to any 
structured cohomology theory whose coefficients are known. We will explicitly emphasize three 
main examples. The first two are to known 
theories, namely ordinary differential cohomology and differential K-theory. We take this opportunity to 
explicitly develop the third theory, which is differential Morava K-theory and then apply our smooth 
AHSS construction to it. 

\subsection{Ordinary differential cohomology theory}
\label{Sec Coh}

We begin by recovering the usual hypercohomology spectral sequence for the Deligne complex 
(see \cite{Bry}, \cite{EV} Appendix) using our methods. 
We will first look at manifolds, then products of these, and then more generally to smooth 
fiber bundles. 

Let us consider the smooth spectrum ${\rm diff}(\Sigma^nH\ZZ,{\rm ch})$ representing differential 
cohomology in degree $n$. We would like to see what our smooth AHSS gives in this case. We recall that 
${\rm diff}(\Sigma^nH\ZZ,{\rm ch})$ is represented  by Deligne cohomology of the sheaf of chain 
complexes $\ZZ^{\infty}_{\cal D}(n)$ via the Eilenberg-MacLane functor 
$
H:Sh_{\infty}(\cartsp; \mathscr{C}{\rm h})\to Sh_{\infty}(\cartsp; \mathscr{S}{\mathrm p})
$
(expressions \eqref{Smooth EM}).  
It follows from the general properties of this functor that the homotopy groups are given by
$$
\pi_{k}{\rm diff}(\Sigma^nH\ZZ,{\rm ch})\simeq H_k\ZZ^{\infty}_{\cal D}(n)\;.
$$
In this case we have the immediate corollary to Theorem \ref{MV thm}. 

\begin{corollary}
The spectral sequence for Deligne cohomology takes the form
$$
E^{p,q}_2=H^p(X;H_{-q}\ZZ^{\infty}_{\cal D}(n))\ \Rightarrow H^{p+q}(X;\ZZ^{\infty}_{\cal D}(n))\;,
$$
which is essentially the hypercohomology spectral sequence for the Deligne complex, but shifted as a fourth quadrant spectral sequence. 
\end{corollary} 
For the sake of completeness, we work out this spectral sequence and recover the differential cohomology diamond 
\eqref{spectra diamond}
from the sequence. This will help to illustrate how the general spectral sequence behaves and how it can be used to calculate general differential cohomology groups.

Now over the site of Cartesian spaces, the Poincar\'e Lemma implies that we have an isomorphism of 
presheaves $d:\Omega^{n-1}/{\rm im(d)}\overset{\simeq}{\to} \Omega^n_{\rm cl}$. Since 
$\Omega^n_{\rm cl}$ is a sheaf over the site of smooth manifolds, the gluing condition allows 
us to calculate the relevant terms on the $E_2$-page of the spectral sequence:
\begin{center}
\begin{tikzpicture}
  \matrix (m) [matrix of math nodes,
    nodes in empty cells,nodes={minimum width=5ex,
    minimum height=4ex,outer sep=-5pt},
    column sep=1ex,row sep=1ex]{
                &      &     &     & 
                \\         
           {\tiny 1} &     &    &   &
            \\
             {\tiny 0}     &   \Omega^{n}_{\rm cl}(M) &   &   & 
             \\   
             \vdots & & & 0 &
             \\       
        {\tiny -(n-2)}\ \ \ \ \ \ \ \ \ & & 0 & & &
             \\
       {\tiny  -(n-1)}\ \ \ \ \ \ \ \ \ & &  &  &  H^{n-1}(M;U(1))   &   & 
                \\
          {\small   -n}\ \ \ \ \ \ &  &    &   &  &  & 0 &
               \\
                &      &     &     & \\};
                
  \draw[-stealth] (m-3-2.south east) -- node[above]{\small $d_2$}(m-4-4.north west);
   \draw[-stealth] (m-6-5.south east) -- node[above] {\small $d_2$}  (m-7-7.north west);
      \draw[-stealth] (m-5-3.south east) -- node[above] {\small $d_2$}  (m-6-5.north west);

\draw[thick] (m-1-1.east) -- (m-8-1.east) ;
\end{tikzpicture}
\end{center}
The term $H^{n-1}(M;U(1))$ will survive to the $E_{\infty}$-page and we have an isomorphism
$$
H^{n-1}(M;U(1))\simeq \frac{F_{n-1}\widehat{H}^n(M;\ZZ)}{F_{n}\widehat{H}^n(M;\ZZ)}\;.
$$
In fact, it is not hard to see that the definition of the filtration gives $F_{n}\widehat{H}^n(M;\ZZ)\simeq 0$ 
and we have an injection
$$
H^{n-1}(M;U(1))\simeq F_{n-1}\widehat{H}^n(M;\ZZ)\into \widehat{H}^n(M;\ZZ)\;.
$$
On the $E_{n}$-page we get one possibly nonzero differential
$$
d_{n}:\Omega^{n}(M)_{\rm cl}\to H^{n}(M;U(1))\;.
$$

\begin{proposition}\label{differential 1}
The differential $d_{n}$ for the AHSS for Deligne cohomology can be identified with the composition
$$
\Omega_{\rm cl}^n(M)\to H^n_{\rm dR}(M)\overset{\int_{\Delta^n}}{\longrightarrow} H^{n}(M;\RR)\overset{\exp}{\longrightarrow} H^{n}(M; U(1))\;,
$$
 and the kernel is precisely those forms which have integral periods.
 \end{proposition}
 \theproof
 We will unpack the definition of the differential in the AHSS in detail. This in turn will require unpacking
 the connecting homomorphism in the Deligne model of ordinary differential cohomology (see \cite{Bry}).  
Denote by $X_p$ the {\v C}ech filtration, and let 
$$
\partial:{\rm diff}(\Sigma^nH\ZZ,{\rm ch})^{q}(X_p)\to {\rm diff}(\Sigma^nH\ZZ,{\rm ch})^{q+1}(X_{p+1}/X_p)
$$ 
denote the connecting homomorphism in the long exact sequence associated to the cofiber sequence $X_p\into X_{p+1}\onto X_{p+1}/X_p$ in the usual way. In what follows, we will denote {\v C}ech-Deligne cochains on the $p$-th level of the filtration $X_p$ as a $p$-tuple
$$(z_0,z_1,\hdots, z_p)\in \widehat{C}^q(X_p)\;,$$
where $z_i$ is a $(q-i)$-form defined on $i$-fold intersections.

Now, by definition, $d_n:E^{0,0}_n\to E^{n,0}_n$ is given by $d_n=\partial (j^*)^{-1}$, where $(j^*)^{-1}$ denotes a choice of element in the preimage of the restriction $j^*$ induced by \footnote{Note that the differential only takes this form at the $(0,0)$-entry. In general, the differential formed from the $n$-th derived couple will be more complicated} $j:X_0\into X_{n-1}$. Since we have $d_k=0$ for $k<n$, the differential 
$d_{n}$ is defined on all elements $z\in \Omega^n_{\rm cl}(M)$. Let $g_0$ be a locally defined $(n-1)$-form trivializing $z$. Then we can choose $(j^*)^{-1}z$ to be the {\v C}ech-Deligne cocycle
 \(\label{cech deligne}
 (j^*)^{-1}z=\underbrace{(g_0,g_1,g_2,\hdots, g_{n-2})}_{n-1}\in \widehat{C}^{0}(X_{n-1})\;,
 \)
 where each $g_k$ is a $(n-k-1)$-form satisfying the cocycle condition $\delta(g_k)=(-1)^{k}dg_{k+1}$.
  To see where the boundary map takes this element, let $y$ be a {\v C}ech-Deligne cochain given by 
 $$
 y=\underbrace{\Big(g_0,g_1,g_2,\hdots,g_{n-2}, {\rm exp}(2\pi i g_{n-1})\Big)\Big)}_{n}\in \widehat{C}^0(X_{n})\;,$$
where $g_{n-1}$ is any smooth $\RR$-valued function satisfying 
\footnote{Note that this cocycle condition is necessary for $y$ to be an lift of $(j^*)^{-1}z$ to the $n$-level of the filtration} $d(g_{n-1})=(-1)^{n-1}\delta(g_{n-2})$. Now $y$ is not {\v C}ech-Deligne closed in general since
$$
 Dy = (d+(-1)^{n-1}\delta) y = (0,0,\hdots, {\rm exp}((-1)^{n-1}2\pi i\cdot \delta(g_{n-1})))
$$
and $g_{n-1}$ may not satisfy the cocycle condition $\delta(g^{n-1})=0$. However, by the 
{\v C}ech-de Rham isomorphism (see for example \cite{BT}), this element in the {\v C}ech-de Rham double complex is isomorphic to an $\RR$-valued {\v C}ech cocycle on $n$-fold intersections. Explicitly, there is a constant $\RR$-valued cocycle $r_n$ such that $\delta(g^{n-1})=r_n$. It follows from the {\v C}ech-singular isomorphism and the singular-de Rham isomorphism that the class of $r_n$ can be represented by the singular cocycle given by the pairing $\int_{\sigma}z$ for any cycle $\sigma$ in $M$. Since the class $\int_{\sigma}z$ was just an unraveling of the boundary $\partial((j^*)^{-1}z)$, we have proved the claim.
  \endofproof

In the next section, we will need to make use of a differential refinement of the Chern character. To this end, we briefly discuss differential cohomology with rational coefficients $\widehat{H}^n(-;\QQ)$. These groups are obtained via the differential function spectra ${\rm diff}(\Sigma^nH\QQ,{\rm ch})$ which fits into the homotopy cartesian square
$$
\xymatrix{
{\rm diff}(H\QQ,{\rm ch})\ar[r]\ar[d] & H(\tau_{\leq 0}\Omega^*[n])\ar[d]
\\
\Sigma^nH\underline{\QQ}\ar[r] & H(\Omega^*[n])
\;.
}
$$
As a consequence of Proposition \ref{diff-cohomology-cal}, the cohomology groups with values in this spectrum are calculated as
$$
{\rm diff}(\Sigma^nH\QQ,{\rm ch})^q(M)=\left\{\begin{array}{cc}
H^{n+q}(M) & q>0,
\\
\\
\widehat{H}^{n}(M;\QQ) & q=0
\\
\\
H^{n-1+q}(M;\RR/\QQ)& q<0\;.
\end{array}\right.
$$
The explicit calculation of the differential in Proposition \ref{differential 1} can be easily modified to get the following. \footnote{The exact argument in the proof of proposition \ref{differential 1} applies, with $\RR/\QQ$ in place of $\RR/\ZZ\simeq U(1)$.}
\begin{proposition}\label{differential rational}
The differential $d_n$ on the $E_n$-page for the AHSS spectral sequence for ${\rm diff}(\Sigma^nH\QQ,{\rm ch})$ is given by 
$$
\Omega_{\rm cl}^n(M)\to H^n_{\rm dR}(M)\overset{\int_{\Delta^n}}{\longrightarrow} H^{n}(M;\RR)\longrightarrow H^{n}(M; \QQ/\ZZ)\;,
$$
 and the kernel is precisely those forms which have rational periods.
\end{proposition}

We will make use of this result when we discuss the differentials in smooth K-theory
in the next section. 
For now, from Proposition \ref{differential 1}, we  immediately get the following characterization of 
 closed forms with integral periods and forms with rational periods using our smooth AHSS.
 \begin{corollary} {\bf (i)}
 The group of closed forms with integral periods on a manifold $M$ is given by 
 $$\Omega^{n}_{\rm cl, \ZZ}(M)\simeq \frac{\widehat{H}^n(M;\ZZ)}{F_1\widehat{H}^{n}(M;\ZZ)}\;.$$
 \item {\bf (ii)} The group of closed forms with rational periods on a manifold $M$ is given by 
 $$\Omega^{n}_{\rm cl, \QQ}(M)\simeq \frac{\widehat{H}^n(M;\QQ)}{F_1\widehat{H}^{n}(M;\QQ)}\;.$$
\end{corollary}

\subsection{Differential $K$-theory}
\label{Sec K}

In this section we examine the smooth AHSS for the differential function spectrum ${\rm diff}(K,{\rm ch})$, 
corresponding to complex $K$-theory. Proposition \ref{diff-cohomology-cal} allows us to calculate 
the cohomology groups on a paracompact manifold $M$ as (see \cite{Lo} \cite{BS} \cite{SSu} \cite{FL})
\(\label{K diff groups}
{\rm diff}(K,{\rm ch})^q(M)=\left\{\begin{array}{ccc}
K^q(M), && q>0,
\\
\\
\widehat{K}^0(M), && q=0,
\\
\\
K^q_{U(1)}(M), && q<0.
\end{array}\right.
\)
Now both groups $K$ and $K_{U(1)}$ are periodic. Indeed, $K_{U(1)}(M)$ fits into an exact sequence
$$
\hdots \to K^{-1}(M)\otimes \RR \to K^{-1}_{U(1)}(M) \to K(M) \to K(M)\otimes \RR \to \hdots\;.
$$
Consequently,  the periodicity of both integral and rational $K$-theory, along with an application of the 
Five Lemma, imply that $K_{U(1)}$ is 2-periodic. In particular, we have 
$$
K^{2q}_{U(1)}(\ast)\simeq U(1) \quad {\rm and} \quad
K^{2q+1}_{U(1)}(\ast)\simeq 0\;, \quad q\in \ZZ\;.
$$
Given the correspondence \eqref{K diff groups}, we see that for a contractible open set $U$, we have 
an isomorphism
$$
{\rm diff}(K,{\rm ch})^{2q+1}(U)\simeq K^{2q}_{U(1)}(\ast)\simeq U(1)
$$
for $q<0$. For degree $0$, the differential cohomology diamond in this case takes the form
$$
\xymatrix @C=5pt @!C{
 & \underset{2k-1}{\prod}\Omega^{2k-1}/{\rm im}(d)\ar[]+<2ex,-2ex>;[rd]^-{a}\ar[rr]^-{d} & &  
 \underset{2k}{\prod}\Omega^{2k}_{\rm cl}\ar[rd] &  
\\
K_{\RR}^{-1} \ar[ru]\ar[rd] & & {\widehat{K}^0}\ar[rd]^-{I} \ar[ru]^-{R}& & K_{\RR}^0\;.
\\
&K_{U(1)}^{-1} \ar[ru]\ar[rr]^-{\beta_{K}} & & K^0 \ar[ru]^{\rm ch} &
}
$$
This implies that for a contractible open set $U$, differential K-theory 
$\widehat{K}^0(U)$ fits into the short exact sequence
$$
0\to \prod_{2k-1}\Omega^{2k-1}/{\rm im}(d)(U)\to \widehat{K}^0(U) \to \ZZ \to 0\;.
$$
Hence, over the site of Cartesian spaces, we have a naturally split short exact sequence of presheaves
$$
0\to \underset{2k-1}{\prod} \Omega^{2k-1}/{\rm im}(d)\to \widehat{K}^0 \to \underline{\ZZ} \to 0\;.
$$
Over that site, the presheaf on the left hand side is actually a sheaf and is naturally isomorphic 
(by Poincar\'e lemma) to the sheaf $\prod_{2k}\Omega^{2k}_{\rm cl}$. We therefore make the identification 
\(
\label{K hat omega}
\widehat{K}^0 \simeq \prod_{2k}\Omega^{2k}_{\rm cl}\oplus \underline{\ZZ}\;.
\)
\begin{remark}
It is important to note that the  identification \eqref{K hat omega} is only true on the site of 
\emph{Cartesian spaces}, which is to say that it holds only locally. On the site of 
smooth manifolds, this is of course not the case.
\end{remark}

Next, since both $\Omega^{2k}_{\rm cl}$ and $\underline{\ZZ}$ are sheaves on the site of smooth manifolds, we can identify the degree $0$ {\v C}ech cohomology with these coefficients with the value of this sheaf on $M$. Isolating the terms on the $E_2$-page which converge to $\widehat{K}^0(M)$, we get
$$
\begin{tikzpicture}
  \matrix (m) [matrix of math nodes,
    nodes in empty cells,nodes={minimum width=5ex,
    minimum height=4ex,outer sep=-5pt}, 
    column sep=1ex,row sep=1ex]{
                &      &     &     & 
                \\         
            1 &     &   &   &
            \\
             0     &   \prod_{2k}\Omega^{2k}_{\rm cl}(M)\oplus \ZZ  &  &   &
             \\           
               -1\  \   &  & H^1(M;U(1)) &  H^2(M;U(1))     & 
                \\
              -2\ \  &   &    &  0 & 0 & 0
               \\
               -3\ \ &  & & &  & H^4(M;U(1))
               \\
               -4\ \ &   &    &   &  &  & 
               \\
                &      &     &     & \\};
                
  \draw[-stealth] (m-3-2.east) -- node[above]{\small $d_2$} (m-4-4.north west);
   \draw[-stealth] (m-4-3.south east) -- (m-5-5.west);
      \draw[-stealth] (m-5-4.south east) -- (m-6-6.north west);
        \draw[-stealth] (m-4-4.south east) -- (m-5-6.north west);
\draw[thick] (m-1-1.east) -- (m-8-1.east) ;
\end{tikzpicture}
$$
We see that all the differentials are zero except for the map labelled $d_2$ above. On the $E_3$-page we get 
$$
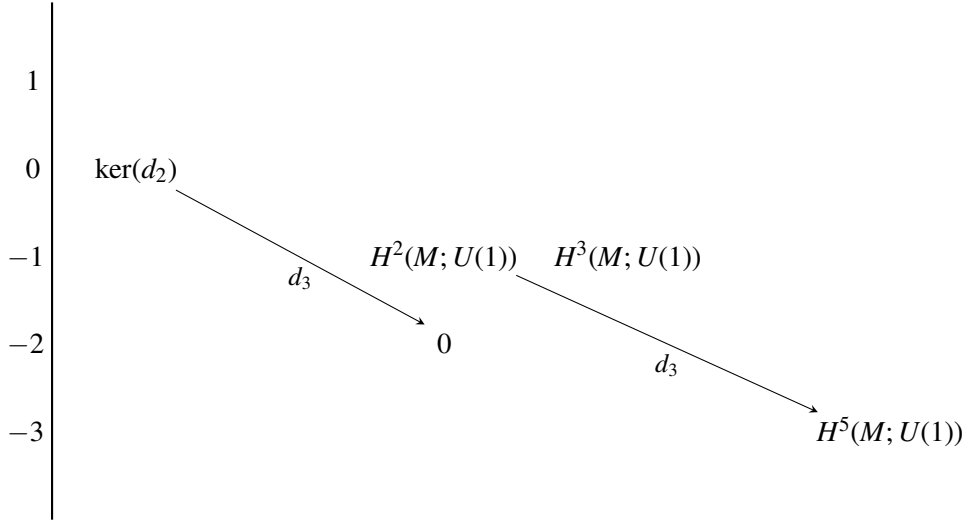
\begin{tikzpicture}
  \matrix (m) [matrix of math nodes,
    nodes in empty cells,nodes={minimum width=5ex,
    minimum height=5ex,outer sep=-5pt}, 
    column sep=1ex,row sep=1ex]{
                &      &     &     & 
                \\         
            1 &     &   &   &
            \\
             0     &   \ker(d_2)  &  &   &
                \\
              -1\ \  &   &    &  & H^2(M;U(1))  & H^3(M;U(1)) 
               \\
               -2\ \ &  &  &  & 0 
               \\
               -3\ \ &   &    &   &  &  & & H^5(M;U(1))
               \\
                 &      &     &     & \\};
                 
             \draw[-stealth] (m-3-2.south east) --node[pos=.5,below]{\small $d_3$} (m-5-5.north west);      
   \draw[-stealth] (m-4-5.south east) --node[pos=.5,below]{\small $d_3$} (m-6-8.north west);
\draw[thick] (m-1-1.east) -- (m-7-1.east) ;
\end{tikzpicture}
$$
The higher pages will fall into cases depending on the parity. We observe that for each even page $E_{2m}$, 
there is one non-zero differential given by $d_{2m}$. For the odd pages the differentials are given by an
odd-degree $U(1)$-cohomology operation.

Note that, in the diagrams, we are interested in the case $p+q=0$, corresponding to diagonal entries. 
Now $p \geq 0$, as the {\v C}ech filtrations are of non-negative degrees, which implies that $q \leq 0$. 
Hence the entries go down the diagonal. Our first goal will be to identify the even differentials $d_{2m}$. 
In order to do this, let us recall that there is a \emph{differential Chern character} map (see \cite{Bun} \cite{Urs}) which is stably given by a morphism of smooth spectra
$$
\widehat{{\rm ch}}:{\rm diff}(K,{\rm ch})\to \prod_{2k} {\rm diff}(\Sigma^{2k}H\QQ,{\rm ch})\;.
$$
Post-composing this map with the projection ${\rm pr}_{2m}$ onto the $2m$-component gives a map of 
smooth spectra
$$
{\rm pr}_{2m}\widehat{{\rm ch}}:{\rm diff}(K,{\rm ch})\to {\rm diff}(\Sigma^{2m}H\QQ,{\rm ch})\;.
$$
Using this map, we can prove the following analogue of Proposition \ref{differential rational}.
\begin{proposition}
The group of permanent cycles in bidegree $(0,0)$ in the AHSS for ${\rm diff}(K,{\rm ch})$ is a subgroup of even degree closed forms with rational periods. That is, we have
$$E^{0,0}_{\infty}\subset \prod_k\Omega^{2k}_{\rm cl,\QQ}(M)\oplus \ZZ\;.$$
\label{even differentials}
\end{proposition}
\theproof
We prove by induction on the even pages of the spectral sequence that, for all $n$, $E^{0,0}_{2n}$ must be a subgroup of 
 \footnote{The differential is $0$ for the odd pages, and so no generality is lost by restricting to the even pages.}
$$
\prod_{2k\leq 2n}\Omega^{2k}_{\rm cl,\QQ}(M)\oplus \prod_{2k>2n}\Omega^{2k}_{\rm cl}(M)\oplus \ZZ\;.
$$
 For the base case, observe that the map 
${\rm pr}_{2}\widehat{{\rm ch}}$ induces a rank 1 morphism of AHSS's and therefore commutes with $d_2$.
 It is straightforward to check, using the definitions, that this leads to the commutative diagram
$$
\xymatrix{
{\prod}_{2k}\Omega^{2k}_{\rm cl}(M)\oplus\ZZ \ar[r]^-{{\rm pr}_{2}}\ar[d]^{d_2} 
& \Omega^2_{\rm cl}(M)\ar[d]^{d^{\prime}_2}
\\
H^2(M;\RR/\ZZ)\ar[r]^{q} & H^2(M;\RR/\QQ)\;.
}
$$
By Proposition \ref{differential 1}, we see that the kernel of $d_2$ must be a subgroup of $\Omega^2_{\rm cl,\QQ}(M)\oplus \prod_{2k>2}\Omega^{2k}(M)\oplus \ZZ$.

Now suppose the claim is true for $d_{2n}$. Again, we have that ${\rm pr}_{2n+2}\widehat{\rm ch}$ 
commutes with $d_{2n+2}$ and we have a commutative diagram
$$
\xymatrix{
\ker(d_{2n}) \ar[r]^-{{\rm pr}_{2n+2}}\ar[d]^{d_{2n+2}} & 
\Omega^{2n+2}_{\rm cl}(M)\ar[d]^{d^{\prime}_{2n+2}}
\\
H^{2n+2}(M;\RR/\ZZ)\ar[r]^-{q} & H^{2n+2}(M;\RR/\QQ)\;.
}
$$
By the induction hypothesis, 
$$
\ker(d_{2n})\subset \prod_{2k\leq 2n }\Omega^{2k}_{\rm cl,\QQ}(M)\oplus \prod_{2k>n}\Omega^{2k}_{\rm cl}(M)\oplus \ZZ,
$$
and the kernel of $d_{2n+2}$ is as claimed.
\endofproof

We now turn to the first odd differential $d_3$. Recall that $\beta$ and $\tilde{\beta}$ denote the Bockstein homomorphisms corresponding to the sequences
$0\to \ZZ \to \RR \overset{\rm exp}{\longrightarrow} U(1)\to 0$
and 
$0\to \ZZ \to \ZZ \to \ZZ/2\to 0$,
respectively. We still  also denote by $\Gamma_2:H^n(-;\ZZ/2)\to H^n(-,U(1))$ the map induced by the representation of $\ZZ/2$ as the square roots of unity and $\rho_2:\ZZ\to \ZZ/2$ as the mod 2 reduction.

\begin{proposition} [Degree three differential] 
The first odd-degree differential in the AHSS for differential K-theory is given by 
$$
d_3=
\left\{\begin{array}{ccc}
\widehat{Sq}^3:= \Gamma_2 Sq^2\rho_2 \beta, && q<0,
\\
\\
Sq^3_{\ZZ}:= \tilde{\beta} Sq^2 \rho_2, && q>0,
\\
\\
0, && q=0.
\end{array}\right.$$
\label{Prop deg3}
\end{proposition}
\theproof
The case for $q=0$ is obvious. For $q>0$, this follows from the fact that the integration map defines an 
isomorphism $I:{\rm diff}(K,{\rm ch})^q(M)\overset{\simeq}{\longrightarrow} K^q(M)$ for $q>0$. Since the
 differential $d_3$ for the  classical AHSS is given by $Sq^3_{\ZZ}$ and the integration map defines an 
 isomorphism of corresponding first quadrant spectral sequences, the case $q>0$ is settled. 

For $q<0$, Corollary \ref{refinement} implies that the Bockstein $\beta$ commutes with the differentials 
on the $E_3$-page. We therefore have
\(\label{d3}
\beta d_3=Sq^3_{\ZZ}\beta=\widetilde{\beta} Sq^3 \rho_2\beta \;.
\)
Rephrasing,  we have the commuting diagram
$$
\xymatrix{
H^{n-1}(M; U(1)) \ar[r]^-{d_3} \ar[d]_\beta &  H^{n+3-1}(M; U(1)) \ar[d]^\beta 
\\
H^n(M; \ZZ) \ar[r]^-{Sq^3_{\ZZ}} & H^{n+3}(M; \ZZ)\;.
}
$$
We now claim that $\tilde{\beta}=\beta\circ \Gamma_2$. Indeed, we have a 
morphism of short exact sequences
$$
\xymatrix{
\ZZ\ar[rr]^{\times 2}\ar[d]^{\rm id} && \ZZ \ar[rr]^{\rho_2}\ar[d]^{\times \pi i}\ar[d] && \ZZ/2\ar[d]^{\Gamma_2}
\\
\ZZ \ar[rr]^{\times 2\pi i} && \RR\ar[rr]^{\rm exp} && U(1)\;.
}
$$
This morphism induces a morphism on the associated long exact sequences on cohomology. 
The homotopy commutativity of the resulting diagram, after delooping once to extend to the left,  
$$
\xymatrix{
\Z/2 \ar[d]_{\Gamma_2} \ar[r]^{\widetilde{\beta}} & B \Z \ar@{=}[d] \\
U(1) \ar[r]^\beta & B\Z
}
$$ 
immediately establishes the claim.

Now it follows from expression \eqref{d3} that $d_3-\Gamma_2 Sq^3 \rho_2\beta$ is in the kernel of $\beta$. 
By exactness of the Bockstein, this implies that it must be in the image of the exponential map, 
${\rm exp} :H^*(-;\RR) \to H^*(-;U(1))$. Hence there is an operation 
$\psi: H^*(-; U(1)) \to H^{*+3}(-; \RR)$ such that 
$$
\phi:={\rm exp}\circ \psi= {\rm exp}(\psi)=d_3-\Gamma_2 Sq^2 \rho_2 \beta\;.
$$
Equivalently, we have a factorization 
$$
\xymatrix{
H^*(-; U(1)) \ar[rd]_\psi \ar[rr]^\phi && H^{*+3}(-; U(1))\;. \\
& H^{*+3}(-; \R) \ar[ur]_{\rm exp}
}
$$
Now the group of natural transformations $H^*(-;U(1))\to H^{*+3}(-;\mathbb{R})$ is in bijective correspondence with $H^{*+3}(K(U(1),*);\RR)\cong 0$. Hence, 
$\psi=0$. \footnote{In the published version, it was erroneously claimed that $\hom(H^n(M;\RR),A)\cong \hom(H^n(M;\QQ),A)$ for any abelian group. This is false, e.g. for $M=S^n$, $A=\QQ$, $\hom(\RR,\QQ)$ is a $\QQ$-vector space of uncountable dimension, while $\hom(\QQ,\QQ)\cong \QQ$. The proof has been corrected.} Consequently, ${\rm exp} \circ \psi=0$, so that $\phi=0$. Therefore, indeed we have 
$$
d_3=\Gamma_2 Sq^3 \rho_2 \beta\;.
$$

\vspace{-8mm}
\endofproof

\begin{remark}
The above proposition suggests that these operations
 are related to some sort of \emph{differential} Steenrod squares. Indeed, this is the case, which has been
  investigated by the authors in \cite{GS2}, with $\widehat{Sq}^3$ being one such operation.  
\end{remark}

Now that we have established the algebraic construction, we turn to investigating the
 convergence of the spectral sequence from a geometric point of view.
 In particular, we immediately observe that the only terms in the spectral sequence which 
 contain information about differential forms are at $q=0$. These terms converge to elements 
 in the filtered graded complex (since $q=0$)
$$
\widehat{K}(M)/F_1\widehat{K}(M)\;.
$$
Since the filtration is given by the {\v C}ech-type filtration on $M$, we see that this quotient contains elements which have nontrivial data on all open sets, intersections, and higher intersections. For the degrees $q<0$, the filtration quotients
$$
F_p\widehat{K}(M)/F_{p+1}\widehat{K}(M)
$$
have trivial data below $p$-intersections.

In fact, it is not too surprising that this occurs. There is a geometric model for reduced $\widehat{K}^0$ 
which is given by the moduli stack $\coprod_{n\in \NN}\BB U(n)_{\rm conn}$ of unitary bundles vector bundles, equipped 
with Hermetian connection. 
More precisely, let ${\rm Vect}_\nabla$ be the moduli stack of complex Hermetian vector bundles with Hermetian connections. 
It was shown in \cite{BNV} that, after taking the Grothendieck group completion, there is a surjection given by the cycle map
$$
{\rm cycl}:{\rm Gr}(\pi_0{\rm Vect}_\nabla(M))\to \widehat{K}^0(M)\;,
$$
which, in our construction, is equivalent to
$$
{\rm cycl}:{\rm Gr}\Big(\pi_0\map\Big(M,\coprod_{n\in \NN}\BB U(n)_{\rm conn}\Big)\Big)\to \widehat{K}^0(M)\;.
$$
Now the stack $\BB U(n)_{\rm conn}$ can be identified with the moduli stack obtained by taking the 
nerve of the action groupoid $
C^{\infty}(-,U(n))//\Omega^1(-;\mathfrak{u}(n))\;,$
with the action given by gauge transformations, where $\frak{u}$ is the Lie algebra of the unitary group. 
Let $\{U_{\alpha}\}$ be a good open cover of $M$. Then a map $
M\to \coprod_{n\in \NN}\BB U(n)_{\rm conn}
$
is given by the  following data:
\begin{list}{$\circ$}{}  
\item A choice of smooth $U(n)$-valued function $g_{\alpha\beta}$ on intersections $U_{\alpha}\cap U_{\beta}$.
\item A choice of local connection 1-form ${\cal A}_{\alpha\beta}$ on open sets $U_{\alpha}$\;.
\end{list}
This is precisely the data needed to define a unitary vector bundle on $M$. 
\begin{remark}
More relevant 
to our needs though, is the fact that the effects of the filtration become transparent 
when taking the completion of $\coprod_{n\in \NN}\BB U(n)_{\rm conn}$ as a model for $\widehat{K}^0$.  We now see that the $q=0$ terms 
converge to terms which involve the data of the connection, while the $q<0$ terms contain data about bundles
 with trivializable connections (in particular, flat connections).
\end{remark}

\paragraph{Differential $K^1$-theory.}
We now consider odd differential K-theory $K^1$. 
In this case the representing spectrum is the unitary group $U$ itself. Viewing this
as a classifying space we can write $U=B \Omega U$.  Of course we are interested in 
 the corresponding stacks.
 Unfortunately, we do not have the analogue of the above group-loop group relation 
 in stacks, i.e. $U_{\rm conn}\not\simeq \mathbf{B}\Omega U_{\rm conn}$. 
Nevertheless, the machinery that we set up will work equally well for differential $K^1$-theory, as far
as the third differential goes, i.e. $d_3=\widehat{Sq}^3$ still. However, the even 
differential are now transgressed in degree by one, so that they are also of 
odd degree. This is expected as the Chern character  in this case is a map to 
 cohomology of odd degree. 

The story for $\widehat{K}^1$ can be worked out similarly, as we indicated above. 
Let us expand on this in more details. 
In the odd case, the differential cohomology diamond takes the form

$$
\xymatrix @C=20pt @!C{
 &\underset{2k}{\prod}\Omega^{2k}/{\rm im}(d)\ar[]+<2ex,-2ex>;[rd]^-{a}\ar[rr]^-{d} & &  
 \underset{2k+1}{\prod}\Omega^{2k+1}_{\rm cl}\ar[]+<3ex,-2ex>;[rd] &  
\\
K_{\RR}^{0} \ar[ru]\ar[rd] & & {\widehat{K}^1}\ar[rd]^-{I} \ar[ru]^-{R}& & K_{\RR}^1
\\
&K_{U(1)}^{0} \ar[ru]\ar[rr]^-{\beta_{K}} & & K^1 \ar[ru]^{\rm ch} &
}
$$
and we get a short exact sequence of presheaves (on the site of Cartesian spaces)
$$
0\to \underline{\ZZ}\to \prod_{2k}\Omega^{2k}/{\rm im}(d)\to \widehat{K}^1\to 0\;.
$$
It is straightforward to show that the map $\underline{\ZZ}\to \prod_{2k}\Omega^{2k}/{\rm im}(d)$ 
is zero. Consequently, we have the isomorphism
$$
\widehat{K}^1\simeq \prod_{2k}\Omega^{2k}/{\rm im}(d)\simeq \prod_{2k+1}\Omega^{2k+1}_{\rm cl}\;.
$$
Using the same type of argument as in the even K-theory $K^0$, we likewise 
get a refinement of the differential of the underlying topological theory. 
More precisely, 
we see that the first nonzero differentials appear on the $E_3$-page is 

$$
\begin{tikzpicture}
  \matrix (m) [matrix of math nodes,
    nodes in empty cells,nodes={minimum width=5ex,
    minimum height=3ex,outer sep=-5pt}, 
    column sep=1ex,row sep=1ex]{
                &      &     &     & 
                \\         
            1 &     &   &   &
            \\
             0     &   \prod_{2k+1}\Omega^{2k+1}_{\rm cl}(M)  &  &   &
             \\           
               -1\  \   &  &  &      & 
                \\
              -2\ \  &   &    &   & H^2(M;U(1)) & H^3(M;U(1))
               \\
               -3\ \ &  & & &  & 
               \\
               -4\ \ &   &    &   &  &  & & H^5(M;U(1))
               \\
                &      &     &     & \\};
                
  \draw[-stealth] (m-3-2.south east) -- node[above]{\small $d_3$} (m-5-6.north west);
        \draw[-stealth] (m-5-5.south) -- node[below]{\small $\widehat{Sq}^3$}(m-7-8.north west);
\draw[thick] (m-1-1.east) -- (m-8-1.east) ;
\end{tikzpicture}
$$
\begin{proposition}
Proposition \ref{Prop deg3} holds for differential $K^1$-theory. That is, the 
degree three differential in $\widehat{K}^1$ is given by the refinement of the 
Steenrod square of dimension three. 
\end{proposition} 

Furthermore, using the same argument as in the proof of  Proposition \ref{even differentials}, 
we see that the permanent cycles in bidegree $(0,0)$ are a subgroup of odd degree forms with rational periods. 

\begin{proposition}
The group of permanent cycles in bidegree $(0,0)$ in the AHSS for ${\rm diff}(\Sigma K,{\rm ch})$ is a subgroup of odd degree closed forms with rational periods. That is, we have
$$E^{0,0}_{\infty}\subset \prod_k\Omega^{2k-1}_{\rm cl,\QQ}(M)\oplus \ZZ\;.$$
\label{prop odd differentials}
\end{proposition}

\begin{example} [Fields in string theory and M-theory] 
In the string theory and M-theory literature one encounters
settings where cohomology classes are compared to K-theory elements, in the
sense of asking when a cohomology class arises from or `lift to' a K-theory class.
This involves, in a sense, a physical modelling of the process of building the 
AHSS. One such obstruction is $Sq^3$, viewed as the first nontrivial differential 
$d_3$ in K-theory, so that the condition $Sq^3 x=0$ on a cohomology class $x$ amounts to 
saying that the class lift to K-theory. This is desirable in the study of the partition 
function of the fields in type IIA string theory (see \cite{DMW} \cite{KS1}). 
On the other hand, it is desirable to have differential refinements for physical 
purposes. Therefore, now that we have the differential AHSS at our disposal, 
it is natural to consider expressions such as 
$d_3 (\widehat{x}):=\widehat{Sq}^3 \hat{x}=0$ on the differential cohomology 
class $\widehat{x}$ that refines the topological class $x$. This can be viewed as 
a condition on cohomology with $U(1)$-coefficients (or flat $n$-bundles), in order 
that they lift to flat elements in $\widehat{K}$. 
\footnote{This could end up being stronger in the 
sense that it is a condition for lifting 
differential cohomology classes to differential K-theory, but we will leave that 
for future investigations.} If the degree of the class
$x$ is even then we are in type IIA string theory and we lift to differential $K^0$-theory. 
On the other hand, being in type IIB string theory means the degree of $x$ is odd, and 
we are lifting to differential $K^1$-theory. The new differentials $d_{2m}$ and $d_{2m+1}$
arising from differential forms will correspond to even and odd degree closed differential forms, 
as the particular forms representing the physical fields $F_{2m}$ and $F_{2m+1}$ 
via the Chern character. 
\label{Ex RR}
\end{example}

\begin{example}[D-brane charges]
The charges of D-branes can a priori be taken to be given as a class in 
cohomology $\cQ_H\in H^*(X; \QQ)$. Quantum effects requires
some of these charges to be (up to shifts) to be in integral cohomology. 
However, in order to not discuss isomorphism classes of 
such physical objects but pinning down a particular physical object, one
considers the charges to take values in differential cohomology, with 
Deligne cohomology being one such presentation, $\cQ_{\hat{H}} \in \widehat{H}^*(X; \Z)$
(see \cite{CJM}).  
On the other hand, careful analysis reveals that the charges 
take values in K-theory rather than in cohomology $\cQ_K \in K^i(X)$, for $i=0, 1$ for
type IIB/IIA (see \cite{MM} \cite{FW} \cite{BMRS}). Such a class exists if 
the cohomology charge satisfies $Sq^3 \cQ_H=0$. 
Again at this stage adding in the 
geometry requires the charges to take values in 
differential K-theory $\cQ_{\widehat{K}} \in \widehat{K}^i(X)$. 
Our construction now allows for a characterization of when
charges in Deligne cohomology lift to charges in differential K-theory,namely 
when they are annihilated by the third differential in the smooth 
AHSS, i.e. when $\widehat{Sq}^3 \cQ_{\hat{H}}=0$. 
\end{example}

\subsection{Differential Morava K-theory} 
\label{Sec Kn}

There are various interesting generalized cohomology theories that descend from 
complex cobordism, among which are Morava K-theory and Morava E-theory. 
Such theories can be defined using their coefficient rings, which in general 
are polynomials over finite or $p$-adic fields on generators whose dimension 
depends on the chromatic level and the prime $p$. As such, these kind of theories
do not lend themselves directly to immediate geometric interpretation, in contrast to 
the case of K-theory, which can be formulated via stable isomorphism classes of 
 vector bundles. 
 
 However, recent work in \cite{LSW} (generalizing some aspects of 
 \cite{BDR}) seems to give hope in that direction.
 Nevertheless,  just because an entity is defined over a finite field does not automatically
 make it ineligible for differential refinement. In fact, recently \cite{GS2} we have demonstrated
 this for the case of Steenrod cohomology operations, which are a priori $\Z/p$-valued
 operations. The main point there was that as long as these admit integral lifts
 then they do have  a chance at a differential refinement. What we will seek here
 is something analogous: integral refinements of such generalized cohomology theories.

We will consider the integral Morava K-theory $\widetilde{K}(n)$, highlighted in \cite{KS1} \cite{S1} \cite{SW}. Morava K-theory $K(n)$ is the mod $p$ reduction of an integral (or $p$-adic) lift $\tK(n)$ with coefficient 
ring $\tK(n)_*=\Z_p[v_n, v_n^{-1}]$.
This theory more closely resembles complex K-theory than is the case for the mod 
$p$ versions (for $n=1$, it is the $p$-completion of K-theory). The integral theory is much more suited to applications in physics \cite{KS1} \cite{S1}  \cite{Buh} \cite{SW}.

The Atiyah-Hirzebruch spectral sequence for Morava K-theory has been studied by Yagita  in \cite{Y} 
(see also \cite{KS1}).   There is a spectral sequence converging to $K(n)^*(X)$ with $E_2$-term 
$E_2^{p, q} = H^p(X, K(n)^q)$. 
While this can be done for any prime, we will focus on the prime $2$. 
 In this case, the first possibly nontrivial differential is $d_{2^{n+1}-1}$; this is given by \cite{Y}
$$
d_{2^{n+1}-1}(xv_n^k) = Q_n(x) v_n^{k-1}\;.
$$
Here $Q_n$ is the $n$th Milnor primitive at the prime $2$, which may be defined inductively as
$Q_0 = Sq^1$, the Bockstein operation, and 
$Q_{j+1}=Sq^{2^j}Q_j - Q_j Sq^{2^j}$, where 
$Sq^j: H^n(X; \Z_2) \to H^{n+j}(X;\Z_2)$ 
is the $j$-th Steenrod square.
These operations are derivations
$$
Q_j(xy)=Q_j(x)y + (-1)^{|x|}x Q_j(y)\;.
$$
The signs are of course irrelevant at $p=2$, but will become important in the integral version.
Extensive discussion of the mod $p$ Steenrod algebra in terms of
 these operations is given in \cite{Tam}.

The integral theory is also computable via an AHSS, 
which   can be deduced from \cite{KS1}  \cite{SW}.
There is an AHSS converging to $\tK(n)^*(X)$ with 
$E_2^{p, q} = H^p(X, \tK(n)^q)$.  The first 
possibly nontrivial differential is $d_{2^{n+1}-1}$; this is given by
$$
d_{2^{n+1}-1}(xv_n^k) = \tQ_n(x) v_n^{k-1}\;.
$$
Here $\tQ_k: H^*(X; \Z) \to H^{*+2^{k+1}-1}(X; \Z)$ is an integral 
cohomology operation lifting the Milnor primitive $Q_k$.

In order to consider differential refinement of Morava K-theory, we need geometric information
encoded in differential forms, hence rational information. 
The rationalization of Morava K-theory $\widetilde{K}(n)$, like any reasonable 
spectrum exists and can be thought of as localization at $\widetilde{K}(0) =H\QQ$. See \cite{Bo} \cite{Ra}. 
We can in the same way localize at $\R$. 
More precisely, the localized theory is given by
$$
\widetilde{K}_{\RR}(n)=\widetilde{K}(n)\wedge M\RR\;,
$$
where $M\RR$ is an Eilenberg-Moore spectrum. We have an equivalence
$$
\widetilde{K}_{\RR}(n)\simeq H\left(\Z[v_n, v_n^{-1}]\otimes \RR\right)
$$
and a Chern character map
$$
{\rm ch}:\widetilde{K}(n)\to  H\left(\Z[v_n, v_n^{-1}]\otimes \Omega^*\right)\;.
$$
Thus we can form the differential function spectrum ${\rm diff}(\widetilde{K}(n),{\rm ch})$ and we can form the associated AHSS. 
To see what form the spectral sequence takes, we need to discuss the {\it flat Morava K-theory} $\widetilde{K}_{U(1)}(n)$, defined by the fiber sequence
$$
\widetilde{K}(n)\to \widetilde{K}(n)\wedge M\RR \to \widetilde{K}_{U(1)}(n):=\widetilde{K}(n)\wedge MU(1)\;.
$$
This theory is periodic with period $2(2^{n}-1)$. Indeed, both $\widetilde{K}(n)$ and its rationalization are 
periodic and we have a long exact sequence
$$
\hdots \tK(n)^m(M)\to  (\tK(n)\wedge M\RR)^m(M)\to \tK_{U(1)}^m(n)(M)\to  \tK(n)^{m+1}(M) \to \hdots \;
$$
relating the flat theory to both the rational ind integral theory.
This, in particular, gives the following identification. 
\begin{lemma} 
 The coefficients of flat Morava K-theory are given by
$$
\tK_{U(1)}(n)^m(\ast)\simeq \left\{\begin{array}{ccl} U(1), && m=2(2^{n}-1),
\\ 
\\
0, && \text{otherwise}.
\end{array}\right.
$$
\end{lemma} 

Knowing the coefficients of the flat theory, we can write down the relevant nonzero terms on the  $E_{2(2^n-1)}$-page of the 
corresponding spectral sequence
$$
\hspace{-15mm}
\scalebox{.9}{
\begin{tikzpicture}
  \matrix (m) [matrix of math nodes,ampersand replacement=\&,
    nodes in empty cells,nodes={minimum width=5ex,
    minimum height=5ex,outer sep=-5pt},
    column sep=.2ex,row sep=1ex]{
                \&      \&     \&     \& 
                \\         
            1 \&     \&   \&   \&
            \\
             0     \&   \prod_{k}\Omega^{k2(2^n-1)}_{\rm cl}(M)\oplus \ZZ  \&  \&   \&
                \\
                \vdots \&
                \\
              -2^{n+1}+3\ \ \ \ \ \ \ \ \ \ \ \ \&   \&    \&  H^{2^{n+1}-3}(M;U(1))  \& H^{2(2^{n}-1)}(M;U(1)) 
               \\
               \vdots  \&  \& \& \& 
               \\
               -2^{n+2}+5 \ \ \ \ \ \ \ \ \ \ \ \ \&   \&    \&   \&  \&  \& H^{4(2^{n}-1)}(M;\RR/\ZZ)
               \\
                \&      \&     \&     \& \\};
                
                \draw[-stealth] (m-3-2.south east)--node[auto]{\small $d_{2(2^n-1)}$} (m-5-5.north west);
                \draw[thick] (m-1-1.east) -- (m-8-1.east) ;
\end{tikzpicture}}
$$
and the only nonzero differential is given by 
$$
d_{2(2^n-1)}:\prod_{k}\Omega^{k2(2^n-1)}_{\rm cl}(M)\oplus \ZZ\to H^{2(2^n-1)}(M;\RR/\ZZ)\;.
$$
Just as in the case for differential $K$-theory (see Propositions  \ref{even differentials} and \ref{prop odd differentials}), we have the following.
\begin{proposition}
The group of permanent cycles in bidegree $(0,0)$ in the AHSS for ${\rm diff}(\widetilde{K}(n),{\rm ch})$ is a subgroup of certain closed forms with rational periods. More precisely, we have
$$E^{0,0}_{\infty}\subset \prod_k\Omega^{2k(2^n-1)}_{\rm cl,\QQ}(M)\oplus \ZZ\;.$$
\end{proposition}

To identify the the {\v C}ech cohomology groups with coefficients in $\widehat{K}(n)^0$, we make the identification (as we did for differential $K$-theory)
$$
\widehat{K}(n)^0\simeq \prod_{k} \Omega^{2k(2^n-1)}_{\rm cl}\oplus \underline{\ZZ}
$$
on the site of Cartesian spaces. Again, using the sheaf condition over smooth manifolds, we have
$$
H^p(M;\widehat{K}(n)^0)\simeq  \prod_{k} \Omega^{2k(2^n-1)}_{\rm cl}(M)\oplus \ZZ\;.
$$
We now consider the differential refinement of the (integrally-lifted) Milnor primitive. 
As before, let $\Gamma_2:H^n(-;\ZZ/2)\to H^n(-; U(1))$ denote the map induced by the representation of $\ZZ/2$ as the square roots of unity and let $\rho_2:\ZZ\to \ZZ/2$ denote the mod 2 reduction. 

\begin{lemma} 
\label{Lem Mil}
The integral Milnor primitive $\tQ_n$ factors through the representation 
$\Gamma_2:\ZZ/2\into U(1)$.
That is, there exists an operation $\widehat{Q}_n$ such that 
$$Q_n\rho_2=\rho_2\tQ_n=\rho_2\beta \Gamma_2 \widehat{Q}_n\;,$$
where $\beta$ is the Bockstein for the exponential sequence. 
\end{lemma} 
\theproof
Recall first that $\rho_2\beta \Gamma_2=\rho_2\tilde{\beta}=Sq^1$, where $\tilde{\beta}$ is the 
Bockstein for the mod 2 reduction sequence. We can therefore rewrite the above equation as
$$
Q_n\rho_2=\rho_2\tQ_n=\rho_2\beta \Gamma_2 \widehat{Q}_n=Sq^1\widehat{Q}_n\;.
$$
and the existence of the class $\widehat{Q}_n$ holds if and only if $Sq^1Q_n\rho_2=0$. On 
the other hand, the existence of the integral lift $\tQ_n$ immediately implies this condition.
\endofproof

Again, let $\beta$ and $\tilde{\beta}$ denote the Bockstein homomorphism 
 corresponding to the sequences
$0\to \ZZ \to \RR\to \RR/\ZZ\to 0$
and 
$0\to \ZZ \to \ZZ \to \ZZ/2\to 0$,
respectively. 
Then the following can be proved in a similar way as we did for Proposition \ref{Prop deg3}
in the case of differential K-theory. 

\begin{proposition} [Odd differentials for Morava AHSS]
The $(2^{n+1}-1)$-differential in the AHSS for differential Morava K-theory is given by 
$$d_{2^{n+1}-1}=\left\{\begin{array}{ccc}
 \Gamma_2 \widehat{Q}_n\rho_2 \beta, && q<0,
\\
\\
\tQ_n, && q>0,
\\
\\
0, && q=0.
\end{array}\right.$$
\label{Prop KQ}
\end{proposition}

\begin{remark}[Odd primes] The above discussion has been for the prime 2, that is, we are considering 
integral Morava K-theory as arising from lifting of the $p=2$ Morava K-theory. We can do the same for odd primes, leading to integral Morava K-theory lifted from an odd prime $p$. A similar discussion follows 
and we have an integral lift of the Milnor primitive at odd primes, as in Lemma \ref{Lem Mil}. The 
differentials will be again given by these refinement of the Milnor primitive, i.e. Proposition \ref{Prop KQ} 
holds except that the primitives are defined using the Steenrod reduced power operations $P^j$. 
Precisely, $Q_0$ is the Bockstein homomorphism associated to reduction mod $p$ sequence, and 
inductively 
$
Q_{i+1}=P^{p^i}Q_i - Q_iP^{p^i}
$.
The operations $P^j$ have been differentially refined in \cite{GS2}. Hence the refinement of the 
Milnor primitives at odd primes will also follow. Then
the $(p^{n+1}-1)$-differential in the AHSS for differential Morava K-theory is given by 
$$d_{p^{n+1}-1}=\left\{\begin{array}{ccc}
 \Gamma_p \widehat{Q}_n\rho_p \beta, && q<0,
\\
\\
\tQ_n, && q>0,
\\
\\
0, && q=0.
\end{array}\right.$$
\end{remark}

\begin{example}[Lifting fields to differential Morava K-theory]
We will build on Example  \ref{Ex RR} and aim to lift the cohomology 
classes beyond K-theory. 
In particular, for $x=\lambda=\tfrac{1}{2}p_1$ the first Spin characteristic class,
we have $\widehat{x}=\hat{\lambda}$ the differential refinements of 
$\lambda$ \cite{SSS3} \cite{Cech} (which can be viewed as a lifted Wu 
class \cite{HS}) we would have $\widehat{Sq}^3 \hat{\lambda}=0$. 
This condition in differential cohomology can be viewed as a refinement of the 
condition $W_7=Sq^3 \lambda=0$ leading to orientation with respect to 
integral Morava K(2)-theory (lifted from the prime $p=2$) as shown in 
\cite{KS1} and elaborated further in \cite{Buh}. From the structure of the smooth 
AHSS in relation to the classical AHSS, one can extend various results to 
the differential case. For instance, one can generalize the statement in \cite{KS1}
on orientation to state that: {\rm an oriented smooth 10-dimensional manifold 
is oriented with respect to differential (integrally lifted from $p=2$) Morava
K(2)-theory $\widehat{K}(2)$ if the class $\hat{W}_7:=\widehat{Sq}^3 \hat{\lambda}=0$.}
The development of this as well as the relation to refinements of characteristic 
classes deserves a separate treatment and will be addressed elsewhere. 
\end{example}

\begin{remark}
{\bf (i)} Note that our construction allows for an AHSS for other spectra beyond the particular 
ones we discussed above. This holds for any spectrum which admits a rationalization,
 whose coefficients are known, and which can be lifted integrally in the sense that we discussed at the 
 beginning of this section. 
\item {\bf (ii)} All the cohomology theories that we used in this paper can be twisted. 
Indeed, the construction in this paper can be generalized to construct
an AHSS for twisted differential spectra \cite{GS4}, in the sense of \cite{BN}.
\end{remark}

\vspace{6mm}
\noindent {\bf \large Acknowledgement}

\noindent The authors would like to thank Ulrich Bunke, Thomas Nikolaus,
and Craig Westerland for interesting discussions at the early stages of this 
project and Urs Schreiber
for very useful comments on the first draft. We are grateful to the anonymous
referee  for a careful reading and for useful suggestions.


\end{document}